%% file: discreteperimeters.tex
\documentclass[10pt]{amsart}
\usepackage{amssymb,amsmath,amsthm,amsfonts}
\usepackage{mathrsfs,dsfont,a4wide}
\usepackage{color}
\usepackage{url}

\usepackage{ifpdf}
\ifpdf
    \usepackage[pdftex]{graphics,graphicx}
    \DeclareGraphicsExtensions{.png,.pdf}
\else
   \usepackage[dvips]{graphics,graphicx}
   \DeclareGraphicsExtensions{.eps}
\fi

\theoremstyle{plain}

\newtheorem{theorem}{Theorem}[section]
\newtheorem{lemma}[theorem]{Lemma}
\newtheorem{proposition}[theorem]{Proposition}

\newtheorem{corollary}[theorem]{Corollary}
\newtheorem{remark}[theorem]{Remark}
\newtheorem{definition}[theorem]{Definition}
\theoremstyle{definition}
\theoremstyle{remark}
\numberwithin{equation}{section}

\newcommand{\hs}{{\mathcal H}}

\newcommand{\dist}{{\textrm{dist}\,}}

\newcommand{\R}{{\mathbb R}}

\newcommand{\N}{{\mathbb N}}

\newcommand{\Z}{\mathbb Z}

\newcommand{\Om}{\Omega}

\newcommand{\weakst}{\stackrel{\ast}{\rightharpoonup}}

\newcommand{\hn}{\hs^{N-1}}

\newcommand{\e}{\varepsilon}
\newcommand{\Persf}{Per_{\Sigma,F}}

\newcommand{\res}{\mathop{\hbox{\vrule height 7pt width .5pt depth 0pt
\vrule height .5pt width 6pt depth 0pt}}\nolimits}

\title
[Continuous limits of discrete perimeters]
{Continuous limits of discrete perimeters}
\author[A. Chambolle]
{Antonin Chambolle}
\address[Antonin Chambolle]{CMAP, \'{E}cole Polytechnique, CNRS 91128, Palaiseau, France}
\email[A. Chambolle]{antonin.chambolle@polytechnique.fr}
\author[A. Giacomini]
{Alessandro Giacomini}
\address[Alessandro Giacomini]{Dipartimento di Matematica, Facolt\`a di Ingegneria, Universit\`a degli Studi di Brescia, Via Valotti 9, 25133 Brescia, Italy}
\email[A. Giacomini]{alessandro.giacomini@ing.unibs.it}
\author[L. Lussardi]
{Luca Lussardi}
\address[Luca Lussardi]{Dipartimento di Matematica, I Facolt\`a di Ingegneria, Politecnico di Torino, c.so Duca degli Abruzzi 24, 10129 Torino, Italy}
\email[L. Lussardi]{luca.lussardi@polito.it}
\begin{document}
\vskip .2truecm
\begin{abstract}
\small{We consider a class of discrete convex functionals which satisfy
a (generalized) coarea formula, and study their limit in the 
continuum.
\vskip .3truecm
\noindent Keywords : generalized coarea formula, total variation, anisotropic perimeter.
\vskip.1truecm
\noindent 2000 Mathematics Subject Classification:
49Q20, 65K10.}
\end{abstract}
\maketitle

\section{Introduction}
In the past ten years, optimization methods for image processing
task have made a lot of progress, thanks to the development of
combinatorial methods (maximal flow/minimal cut, and other
graph-based optimization methods --- see for
instance~\cite{GreigPorteousSeheult},
and~\cite{NetworkFlows,BoykovKolmogorov2004} and the references
therein). These methods are not new, the idea of representing Ising
energies (i.e., discrete approximations of perimeters) on graphs 
and computing minimum points using maximal flows algorithms dates back
at least to the 70s~\cite{PicardRatliff}. However, the evolution of
computers and development of new algorithms~\cite{BoykovKolmogorov2004},
oriented towards specific applications, have contributed a lot
towards the recent increase in activity in this field. In image
processing, the idea is to regularize ill-posed inverse problems
for finding sets (shapes) or partitions into labels of an image,
by penalizing a discrete variant of their perimeter. We try to
consider, in this paper, the most general energies which can
be tackled by these methods, and even a little bit more:
we consider discrete \textit{submodular} energies, defined
on discrete subsets of a finite
lattice $\mathcal{V}\subset h\Z^N$, $h>0$, for which
it is known that polynomial algorithms do exist (see for
instance \cite{cunningham.85.combinatoria,iwata.00.acm,schrijver.00.jct}).
We will show that, appropriately extended into functions of
general vectors in $\R^\mathcal{V}$ by means of the generalized
coarea formula, these energies are, in fact, convex. This is
already known (although our setting is a bit different, as
well as our proofs which apply to other situations, including functionals
defined in the continuous setting) in discrete optimization,
under the notion of Lov\'asz' extension~\cite{Lovasz}.

We will then study the continuous limit of our energies, as
the discretization step $h$ goes to 0 (and the number of pixels/voxels
in $\mathcal{V}$ to infinity), providing a very simple representation
formula for the limit. In particular, it will be obvious from this
formula that simple approximation procedures only provide
``crystalline'' energies, as already observed for instance
in~\cite{BK-geocut}.

To be more specific, we consider in this paper an ``interaction
potential'' $F:\{0,1\}^\Sigma\to [0,+\infty)$, which is a nonnegative
function of binary vectors of $\{0,1\}^\Sigma$, where $\Sigma\subset \Z^N$
is a finite (small) set of ``neighbors''. We assume, in addition,
that $F$ satisfies the submodularity condition
$$
F(u \wedge v)+F(u \vee v) \leq F(u)+F(v)
$$ 
for any $u,v \in \{0,1\}^\Sigma$,
where $(u \wedge v)_i:=\min \{u_i,v_i\}$ and $(u \vee v)_i:=\max\{u_i,v_i\}$.
Defining, for $x\in \R^N$ and $u$ a real-valued function, the vector
$u[x+h\Sigma]=(u(x+hi)_{i\in \Sigma})\in \R^\Sigma$, we will study
the asymptotic behavior as $h\to 0$ of functionals of the type
\begin{equation}
\label{discrperintr}
J_h(E,\Om):=h^{N-1}\sum_{x \in I^h(\Om)\cap h\Z^N}F(\chi_E[x+h\Sigma]),
\end{equation}
where here, $\Om$ is a bounded open subset of $\R^N$
and $E$ is a discrete subset of
the discrete lattice $\mathcal{V}=h\Z^N\cap \Om$
($E$ is also identified to the union of the cubes $Q^h_x=x+[0,h[^N$,
$x\in E$, and $\chi_E$ is its characteristic function).
The notation $I^h(\Om)$ stands for the points $x$ such
that $x+h\Sigma\subset\Om$, so that the
sum in \eqref{discrperintr} involves only the nodes $x\in h\Z^N$ such that
$x+h\Sigma \subset \Om$.
\par
The functional \eqref{discrperintr} is a sort of {\it nonlocal anisotropic discrete perimeter} of $E$. In fact it penalizes the boundary of $E$ in a nonlocal way, since 
an interface at the boundary with a vertex $x$ interacts with the behavior of $E$ on the cubes with vertices $x+h\Sigma$. The nonlocality vanishes as $h\to 0$ since its radius of action is given by $h\ diam(\Sigma)$. The anisotropy is introduced by the function $F$, which can weight interfaces with various orientations in different ways.
\par
The main result of the paper concerns the asymptotic behavior of the discrete perimeters \eqref{discrperintr} as $h\to 0$ in the variational sense of $\Gamma$-convergence (see Section \ref{prelsec} for the definition) with respect to the $L^1$-topology on the family of discrete sets (that is $L^1$ convergence of characteristic functions).
Under mild assumptions on $F$ and $\Om$, we prove that (see Theorem \ref{main})  the discrete perimeters $\Gamma$-converge to the continuous anisotropic
perimeter which for a sufficiently regular set $E$ (a set with Lipschitz boundary for instance) is given by
\begin{equation}
\label{contper}
J(E,\Om)=\int_{\partial E}F(\nu_E\cdot \Sigma)\,dA,
\end{equation}
where $\nu_E$ is the inner normal at the boundary and $ (\nu_E\cdot \Sigma)=(\nu_E\cdot y)_{y\in \Sigma}$. More precisely the class on which the $\Gamma$-limit is defined is given by the family of sets with finite perimeter 
in $\Om$. As a consequence, for a general set $E$, the boundary involved in the functional \eqref{contper} is the reduced boundary $\partial^*E$, the inner normal $\nu_E$ is intended in a  measure theoretical sense (see Section \ref{prelsec}), and the area measure $dA$ has to be replaced by the $(N-1)$-dimensional Hausdorff measure $\hn$. The function $F$ appearing in \eqref{contper} is the extension to $\R^\Sigma$ of the submodular function $F$ by means of the formula
\begin{equation}
\label{discretecoareaintr}
F(u)=\int_{-\infty}^{+\infty}F(\chi_{\{u>s\}})\,ds,
\end{equation}
where  $\{u>s\}:= \{x\in\Sigma : u(x)>s\}$. Formula \eqref{discretecoareaintr} is a {\it coarea formula} for the function $F$ since it relates the value $F(u)$ to the behavior of $F$ on the ``boundary'' of the level sets $\{u>s\}$ (compare with equation \eqref{coareaclassic} which gives the classical coarea formula for functions of bounded variation).
\par
In view of the result on dicrete perimeters, we obtain a $\Gamma$-convergence result for the functionals $J_h(\cdot,\Om)$ extended to the class of piecewise constant functions $u$ relative to the grid $h\Z^N$. More precisely we consider $u$ of the form
$$
u=\sum_{x\in h\Z^N}a_x \chi_{Q^h_x},
\qquad a_x\in \R
$$
and
\begin{equation}
\label{Jhintrfunct}
J_h(u,\Om)=h^{N-1}\sum_{x \in I^h(\Om)\cap h\Z^N}F(u[x+h\Sigma]).
\end{equation}
As the functional \eqref{discrperintr} could be thought as a discrete perimeter, the functional \eqref{Jhintrfunct} could be considered as a sort of {\it discrete total variation} of the function $u$. Clearly it inherits the nonlocal and anisotropic features of the discrete perimeter.
We show that the $\Gamma$-limit in the $L^1$-topology is given by the anisotropic total variation
\begin{equation}
\label{totvarintr}
J(u,\Om)=\int_\Om F\left( \frac{Du}{|Du|}\cdot \Sigma\right)\,d|Du|,
\end{equation}
where $u$ belongs to the space $BV(\Om)$ of functions of bounded variation (see Section \ref{prelsec}),
$Du/|Du|\in \mathbb{S}^{N-1}$ denoting the Radon-Nikodym derivative of $Du$
with respect to its total variation $|Du|$. This $\Gamma$-convergence result is a simple consequence of the result on sets and of the fact that the functionals $J_h$ satisfy the {\it generalized coarea formula}
\begin{equation}
\label{coareaJhintr}
J_h(u,\Om)=\int_{-\infty}^{+\infty}J_h(\chi_{\{u>s\}},\Om)\,ds,
\end{equation}
so that the behavior of $J_h$ on piecewise constant functions is completely determined by the discrete perimeters for sets. We infer the result from general  properties of functionals on $L^1(\Om)$ that satisfy a coarea formula like \eqref{coareaJhintr}, which we study in Section \ref{coareasec}. This class of functionals, denoted by $GC(\Om)$, was investigated by Visintin \cite{V1,V2} in connection with phase transition problems. As a consequence of our $\Gamma$-convergence result, the discrete total variations \eqref{Jhintrfunct} can be used to approximate Total Variation Minimization procedures in image denoising involving \eqref{totvarintr} (see Corollary \ref{cormain}).
\par
The paper is organized as follows. Section \ref{prelsec} contains the notation employed in the paper, and some basic facts concerning sets with finite perimeters, functions of bounded variation and $\Gamma$-convergence. In Section \ref{coareasec} we consider the class $GC(\Om)$ of functionals on $L^1(\Om)$ which satisfy the generalized coarea formula: In particular we prove that $GC(\Om)$ is closed under $\Gamma$-convergence, and that the limit can be recovered by the behavior on characteristic functions of Borel sets. Section \ref{discrsec} contains
the main $\Gamma$-convergence result formulated for the discrete total variations \eqref{Jhintrfunct}. We exploit the reduction to the class of discrete sets, and
Subsections \ref{secglinf} and \ref{secglsup} contain the proof of the 
two inequalities characterizing $\Gamma$-convergence for study the discrete
perimeters.

\section{Notation and preliminaries}
\label{prelsec}
Let $A$ be an open subset of $\R^N$. We will say that $A$ has a {\it continuous} boundary if $\partial A$ can be covered by finitely many balls $B$ such that, in each ball, $B\cap A$ is the subgraph of a continuous function (after an appropriate change of coordinates). If these functions are Lipschitz continuous, we say that $A$ has a Lipschitz boundary.

\par
For any $p \in [1,+\infty[$ we will denote by $L^p(A)$ the usual space of all $p$-summable functions on $A$, and by $L^\infty(A)$ the space of measurable functions on $A$ which are essentially bounded. Given $u,v \in L^1(A)$, we set
\begin{equation}
\label{minmax}
u \wedge v:=\min\{u,v\}
\qquad\text{and}\qquad
u\vee v:=\max\{u,v\}.
\end{equation}

In the following, we recall some basic facts concerning function of bounded variation and sets with finite perimeter which we need in the following sections, together with some basic definitions and results concerning $\Gamma$-convergence.

\vskip10pt\noindent
{\bf Functions of bounded variation and sets with finite perimeter.} 
For an exhaustive treatment of the subject, we refer the reader to \cite{AFP}.
\par
We say that $u$ has {\it bounded variation} in $A$ and we write $u \in BV(A)$ if $u \in L^1(A)$ and
\begin{equation}
\label{vartot}
|Du|(A)=\sup\left\{\int_A u\, \textrm{\rm div}\varphi\,dx : \varphi \in C^1_c(A)\,,||\varphi||_\infty \leq 1\right\}<+\infty.
\end{equation}
$|Du|(A)$ is referred to as the {\it total variation} of $u$. 
\par
If $E \subseteq A$ is a Borel set, we say that $E$ has {\it finite perimeter} in $A$ if 
$\chi_E \in BV(A)$, and we set
$$
Per(E,A):=|D\chi_E|(A).
$$
$Per(E,A)$ is called the {\it perimeter} of $E$ in $A$. It turns out that
$$
Per(E,A)=\hn(\partial^* E\cap A),\ \textrm{ and }\ 
D\chi_E\ =\ \nu_E\hn\res\partial^* E,
$$
where $\partial^*E$ denotes the {\it reduced boundary} of $E$,
which, up to a $\hn$-negligible set, coincides with the (larger)
set of points $x$ such that there exists a unit vector $\nu_E(x)$ with
$$
\frac{E-x}{\varrho} \to \{y \in \R^N: y\cdot \nu_E(x)>0\}
\qquad
\text{as $\varrho \to 0$ in } L^1_{loc}(\R^N).
$$
The unitary vector $\nu_E(x)$ is usually referred to as the {\it interior normal} to $E$ at $x$.
$\hn$ denotes the $(N-1)$-dimensional Hausdorff measure, which is a generalization to arbitrary sets of the usual $(N-1)$-area measure. The points of $\partial^*E$ are also called {\it regular points} of $\partial E$.
\par
If $u \in BV(A)$, the following {\it coarea formula} holds:
\begin{equation}
\label{coareaclassic}
|Du|(A)=\int_{-\infty}^{+\infty}Per(\{x \in A : u(x)>s\},A)\,ds
=\int_{-\infty}^{+\infty}\mathcal |D\chi_{\{ u>s\}}|(A)\,ds.
\end{equation}
Finally we recall the following compactness result
(which is a variant of Rellich's theorem).
If $A$ is bounded and with Lipschitz boundary, and $(u_n)_{n \in \N}$ is a sequence in $BV(A)$ such that $\|u_n\|_{L^1(A)}+|Du_n|(A)$ is bounded, then there exist a subsequence $(u_{n_k})_{k \in \N}$ and a function 
$u \in BV(A)$ such that
$$
u_{n_k} \to u \qquad\text{ in }L^1(A)
$$ 
and
$$
|Du|(A) \le \liminf_{k \to \infty} |Du_{n_k}|(A).
$$

\vskip10pt\noindent
{\bf $\Gamma$-convergence.} Let us recall the definition and some basic properties
of De Giorgi's {\it $\Gamma$-convergence} in metric spaces.
We refer the reader to \cite{DalMaso,Braides} for an exhaustive
treatment of this subject. Let $(X,d)$ be a metric space. We say that a sequence
$F_n:X\to [-\infty ,+\infty]$ $\Gamma $-converges to $F:X\to [-\infty ,+\infty ]$
(as $n\to \infty$) if for all $u \in X$ we have
\begin{itemize}
\item[{\rm (i)}] ({\it $\Gamma${-}lim inf inequality}) for every
sequence $(u_n)_{n \in \N}$ converging to
$u$ in $X$,
$$
\liminf_{n \to \infty} F_n(u_n)\geq F(u);
$$
\item[{\rm (ii)}] ({\it $\Gamma${-}lim sup inequality})
there exists a sequence
$(u_n)_{n \in \N}$ converging to $u$ in $X$, such that
$$
\limsup_{n \to \infty} F_n(u_n)\leq F(u).
$$
\end{itemize}
The function $F$ is called the $\Gamma${-}limit of
$(F_n)_{n \in \N}$ (with respect to $d$). Given a family $(F_h)_{h>0}$ of functionals on $X$, we say that $F_h$ $\Gamma$-converges to $F$ as $h \to 0$ if for every sequence $h_n \to 0$ we have that $F_{h_n}$ $\Gamma$-converges to $F$ as $n \to \infty$.
\par
$\Gamma${-}convergence is a convergence of
variational type as explained in the
following Proposition.

\begin{proposition}\label{Gamma-conv-prop}
Assume that the family $(F_h)_{h>0}$
$\Gamma${-}converges to $F$ and
that there exists a compact set $K\subseteq X$ 
such that for all $h>0$
$$
\inf\limits _{u\in K}F_h(u)=\inf\limits _{u\in X}
F_h(u).
$$
Then $F$ admits a minimum on $X$, $\inf_{X}F_h \to \min_X F$ as $h \to 0$, and
any limit point of any sequence $(u_h)_{h>0}$ such that
$$
\lim\limits _{h \to 0}\Bigl( F_h(u_h)-
\inf\limits _{u\in X}F_h(u)\Bigr) =0
$$
is a minimizer of $F$.
\end{proposition}

\section{Generalized coarea formula}
\label{coareasec}
In the following, let $\Omega \subset \R^N$ be an open and bounded set.
\begin{definition}
\label{generalizedcoarea}
Let $J \colon L^1(\Om) \to [0,+\infty]$ be a proper functional. We say that $J$ satisfies the {\sl generalized coarea formula} if for every $u \in L^1(\Om)$

\begin{equation}\label{coarea}
J(u)=\int_{-\infty}^{+\infty} J({\chi_{\left\{u>s\right\}}})\,ds,
\end{equation}
with the convention $J(u)=+\infty$ if the map $s \mapsto J({\chi_{\left\{u>s\right\}}})$ is not measurable. We denote by $GC(\Om)$ the class of functionals satisfying \eqref{coarea}.
\end{definition}

The class $GC(\Om)$ has been introduced by Visintin \cite{V1} and investigated, in the discrete case, by Chambolle and Darbon \cite{CD}.
In a slightly different setting, the formula~\eqref{coarea} is a variant
of the extension introduced by Lov\'asz in~\cite{Lovasz} and
well-known in combinatorial and linear optimization.

An example of functional satisfying \eqref{coarea} is given by the total variation \eqref{vartot} in view of the coarea formula \eqref{coareaclassic}.
Other examples are treated in \cite{V2}:
$$
J(u)=\int_{\Om \times \Om}|u(x)-u(y)||x-y|^{-(N+r)}\,dx\,dy, \quad \forall r \in (0,1)
$$
and
$$
J(u)=\int_{\Om \times \R^+}\big({\underset{B_h(x)\cap\Om}{\textrm{\rm ess sup}}\,u}-
{\underset{B_h(x)\cap\Om}{\textrm{\rm ess inf}}\,u}\big)
h^{-(1+r)}\,dx\,dh, \quad \forall r \in (0,1).
$$
The next Proposition contains some elementary consequences of formula (\ref{coarea}).

\begin{proposition}\label{prop-1}
Let $J \in GC(\Om)$. Then for every $u\in L^1(\Om)$ the following facts hold:
\begin{itemize}
\item[({\rm i)}] $J(\lambda u)=\lambda J(u)$ for every $\lambda > 0$;
\item[({\rm ii)}] $J(u+c)=J(u)$ for every $c \in \R$;
\item[({\rm iii)}] $J(c)=0$ for every $c \in \R$.
\end{itemize}
Moreover if $J$ is convex, for every $u,v \in L^1(\Om)$ we have
\begin{itemize}
\item[({\rm iv)}] $J(u\wedge v)+J(u\vee v) \leq J(u)+J(v)$.
\end{itemize}
\end{proposition}
\begin{proof}
Let $\lambda > 0$, $u \in L^1(\Om)$ and $c \in \R$. Then
$$
J(\lambda
u)=\int_{-\infty}^{+\infty}J({\chi_{\left\{u>\frac{s}{\lambda}\right\}}})\,ds=
\lambda
\int_{-\infty}^{+\infty}J({\chi_{\left\{u>t\right\}}})\,dt=\lambda J(u)
$$
and
$$
J(u+c)=\int_{-\infty}^{+\infty}J({\chi_{\left\{u>s-c\right\}}})\,ds=
\int_{-\infty}^{+\infty}J({\chi_{\left\{u>t\right\}}})\,dt=J(u)
$$
so that (i) and (ii) follow.
\par
Let us prove (iii). In view of (ii), it suffices to show that $J(0)=0$. Suppose by contradiction that $J(0)>0$. Then for every $u \in L^\infty(\Om)$ we have
\begin{equation}
\label{contr}
J(u)=\int_{-\infty}^{+\infty}J(\chi_{\{u>s\}})\,ds
\geq\int_{\underset{\Om}{\textrm{\rm ess sup}}\,u}^{+\infty}J(0)\,ds=+\infty.
\end{equation}
By the generalized coarea formula, we deduce that $J(u)=+\infty$ for every $u \in L^1(\Om)$. But this is against the fact that $J$ is proper, so that point (iii) is proved.
\par
Let us show (iv). Since we have
\begin{align*}
J(\chi_{\{u>s\}} \wedge \chi_{\{v>s\}})+J(\chi_{\{u>s\}}\vee \chi_{\{v>s\}})&=
\int_0^2 J( \chi_{\{\chi_{\{u>s\}} \wedge \chi_{\{v>s\} }+ \chi_{\{u>s\}} \vee \chi_{\{v>s\} } >t\}})\,dt\\
&=\int_{-\infty}^{+\infty} J( \chi_{\{\chi_{\{u>s\}} \wedge \chi_{\{v>s\} }+ \chi_{\{u>s\}} \vee \chi_{\{v>s\} } >t\}})\,dt,\\
\end{align*}
by the generalized coarea formula \eqref{coarea} we get
\begin{multline*}
J(\chi_{\{u>s\}} \wedge \chi_{\{v>s\}})+J(\chi_{\{u>s\}}\vee \chi_{\{v>s\}})=
J(\chi_{\{u>s\}} \wedge \chi_{\{v>s\} }+ \chi_{\{u>s\}} \vee \chi_{\{v>s\}})\\
=J(\chi_{\{u>s\}}+ \chi_{\{v>s\} }).
\end{multline*}
Notice that if $J$ is convex, by point (i) we deduce that $J$ is subadditive. Then we obtain
\begin{equation}
\label{basicineq}
J(\chi_{\{u>s\}} \wedge \chi_{\{v>s\}})+J(\chi_{\{u>s\}}\vee \chi_{\{v>s\}})\le
J(\chi_{\{u>s\}})+J(\chi_{\{v>s\}}).
\end{equation}
Observe that for any $s \in \R$ we have $\{u\wedge v>s\}=\{u>s\}\cap \{v>s\}$ and $\{u\vee v>s\}=\{u>s\}\cup \{v>s\}$ so that
$$
\chi_{\{u \wedge v>s\}}=\chi_{\{u>s\}}\wedge \chi_{\{v>s\}}
\qquad\text{and}\qquad
\chi_{\{u\vee v>s\}}=\chi_{\{u>s\}}\vee \chi_{\{v>s\}}.
$$ 
We conclude by \eqref{basicineq}
$$
\begin{aligned}
J(u \wedge v)+J(u \vee v) &=\int_{-\infty}^{+\infty} [J(\chi_{\{u \wedge v>s\}})+
J(\chi_{\{u \vee v>s\}})]\,ds\\
&=\int_{-\infty}^{+\infty}[J(\chi_{\{u>s\}} \wedge \chi_{\{v>s\}})+J(\chi_{\{u>s\}}\vee \chi_{\{v>s\}})]\,ds\\ &
\leq \int_{-\infty}^{+\infty} J(\chi_{\{u>s\}})\,ds+\int_{-\infty}^{+\infty} J(\chi_{\{v>s\}})\,ds
=J(u)+J(v)
\end{aligned}
$$
so that (iv) follows and the proof is complete.
\end{proof}

We will need the following Lemma concerning the approximation of Lebesgue integral by means of Riemann sums.

\begin{lemma}\label{appross}
Let $f \in L^1(\R)$, $t \in ]0,1[$ and let us set
$$
s_n(t):=\frac{1}{n}\sum_{k \in \Z}f\left(\frac{k+t}{n}\right).
$$
Then up to a subsequence we have 
$$
\lim_{n \to \infty}s_n(t)=\int_{-\infty}^{+\infty} f(\tau)\,d\tau
\qquad\text{for a.e.\,\,$t \in ]0,1[$}.
$$ 
\end{lemma}

\begin{proof}
For any $n \in \N$ we easily get
$$
\int_{-\infty}^{+\infty} f(\tau)\,d\tau=\frac{1}{n}\sum_{k \in \Z}\int_0^1f\left(\frac{k+r}{n}\right)\,dr.
$$
Then for $t \in ]0,1[$ we have
\begin{multline*}
\int_0^1 \left|\int_{-\infty}^{+\infty}f(\tau)\,d\tau-s_n(t)\right|\,dt
=\int_0^1 \left|\int_{-\infty}^{+\infty}f(\tau)\,d\tau-\frac{1}{n}\sum_{k \in \Z}f\left(\frac{k+t}{n}\right)\right|\,dt
\\
=
\int_0^1 \left|\frac{1}{n}\sum_{k \in \Z}\int_0^1f\left(\frac{k+r}{n}\right)\,dr-\frac{1}{n}\sum_{k \in \Z}f\left(\frac{k+t}{n}\right)\right|\,dt 
\\
\le
\frac{1}{n}\sum_{k \in \Z} \int_0^1\int_0^1 \left|f\left(\frac{k+r}{n}\right)-f\left(\frac{k+t}{n}\right) \right|\,dr\,dt.
\end{multline*}
But
\begin{multline*}
\frac{1}{n}\sum_{k \in \Z} \int_0^1\int_0^1 \left|f\left(\frac{k+r}{n}\right)-f\left(\frac{k+t}{n}\right) \right|\,dr\,dt
\\
\le 
\frac{1}{n}\sum_{k \in \Z} \int_0^1\int_0^1 \left[\left|f\left(\frac{k+r}{n}\right)-f\left(\frac{k+r+t}{n}\right) \right|
+\left|f\left(\frac{k+r+t}{n}\right)-f\left(\frac{k+t}{n}\right) \right|
\right]\,dr\,dt\\
=\int_0^1 \left |\left|f(\cdot)-f\left(\cdot+\frac{t}{n}\right)\right |\right|_{L^1(\R)}\,dt+
\int_0^1 \left|\left|f(\cdot)-f\left(\cdot+\frac{r}{n}\right)\right|\right|_{L^1(\R)}\,dr.
\end{multline*}
The last terms tend to zero by continuity of the translation operator in $L^1(\R)$.
We conclude that 
$$
s_n \to  \int_{-\infty}^{+\infty} f(\tau)\,d\tau
\qquad\text{ in }L^1(0,1)
$$
so that, up to a subsequence, pointwise almost everywhere convergence follows.
\end{proof}

In view of \eqref{coarea}, functionals in the class $GC(\Om)$ are completely determined by their behavior on characteristic functions of Borel sets contained in $\Om$. The next result  gives a sufficient condition for the convexity of lower semicontinuous functionals in $GC(\Om)$ in terms of the submodularity property (iv) of the previous Proposition only on characteristic functions.

\begin{proposition}\label{prop-conv}
Let $J \in GC(\Om)$ be a lower semicontinuous functional such that
\begin{equation}
\label{sottomodins}
J(\chi_{E \cap E'})+J(\chi_{E \cup E'}) \leq J(\chi_E)+J(\chi_{E'})
\end{equation}
for every pair of Borel sets $E,E'$ in $\Om$. Then $J$ is convex.
\end{proposition}

\begin{proof}
Since by Proposition \ref{prop-1} $J$ is positively one-homogeneous , it is sufficient to show that 
\begin{equation}
\label{sub2}
J(u+v) \leq J(u)+ J(v)
\end{equation}
for any $u,v \in L^1(\Om)$. 
\par
We claim that the following representation formula holds for every function $u$ which is positive, bounded and with integer values:
\begin{equation}
\label{subadd}
J(u)=\min\left\{\sum_{i=1}^m  J(\chi_{E_i}): m \geq 1\,, u=\sum_{i=1}^m
\chi_{E_i}\right\}.
\end{equation}
In order to prove \eqref{sub2}, we can clearly assume that $J(u)<+\infty$ and $J(v)<+\infty$. Hence by \eqref{coarea} the maps $s \mapsto J(\chi_{\{u>s\}})$ and $s \mapsto J(\chi_{\{v>s\}})$ belong to $L^1(\R)$. 
\par
Firstly let us assume $0\leq u \leq 1$ and 
$0 \le v\le 1$. For every $n \in \mathbb{N}$, $n>0$, let us set for $t \in ]0,1[$
$$
u_n:=\frac{1}{n}\sum_{\begin{subarray}{c} k\in\mathbb{N}
\end{subarray}}\chi_{\left\{u>\frac{k+t}{n}\right\}} \quad\text{and}\quad
v_n:=\frac{1}{n}\sum_{\begin{subarray}{c} k\in\mathbb{N}
\end{subarray}}\chi_{\left\{v>\frac{k+t}{n}\right\}}.
$$
By applying Lemma \ref{appross}, we can choose $t\in ]0,1[$ in such a way that
$$
\lim_{n \to \infty}\frac{1}{n}\sum_{\begin{subarray}{c} k\in\mathbb{N}
\end{subarray}} J\left(\chi_{\left\{u>\frac{k+t}{n}\right\}}\right)=\int_0^1
J(\chi_{\{u>s\}})\,ds 
$$
and
$$
\lim_{n \to \infty}\frac{1}{n}\sum_{\begin{subarray}{c} k\in\mathbb{N}
\end{subarray}}J\left(\chi_{\left\{v>\frac{k+t}{n}\right\}}\right)=\int_0^1
J(\chi_{\{v>s\}})\,ds.
$$
By construction $u_n \to u$ and $v_n \to v$ in $L^1(\Om)$. Then by positive homogeneity, and assuming the representation formula \eqref{subadd} holds,
we get
$$
J(u_n+v_n)=J\left(\frac{1}{n}\sum_{\begin{subarray}{c}
k\in\mathbb{N} \end{subarray}}\chi_{\left\{u>\frac{k+t}{n}\right\}}+\chi_{\left\{v>\frac{k+t}{n}
\right\}}\right)\leq \frac{1}{n}\sum_{\begin{subarray}{c} k\in\mathbb{N}
\end{subarray}}
J\left(\chi_{\left\{u>\frac{k+t}{n}\right\}}\right)+ J
\left(\chi_{\left\{v>\frac{k+t}{n} \right\}}\right).
$$
The right-hand side converges by construction to $J(u)+ J(v)$, and thus, by lower
semicontinuity of $J$, we have that \eqref{sub2} follows.
\par
In the case $m \leq u \leq M$ and $m \le v \le M$, one can easily show again
that \eqref{sub2} holds by considering the functions $(u-m)/(M-m)$ and
$(v-m)/(M-m)$, and taking into account the general properties of $J$.
\par
Finally, for $u,v \in L^1(\Om)$ and for $T>0$, let us consider $u_T:=-T
\vee u \wedge T$ and $v_T:=-T \vee v \wedge T$. Since $u_T \to u$ and
$v_T \to v$ in $L^1(\Om)$, by the lower semicontinuity of $J$ we obtain
\begin{multline*}
J(u+v) \leq \liminf_{T \to +\infty}J(u_T+v_T) \le  \liminf_{T \to +\infty}(J(u_T)+J(v_T)) \le \limsup_{T \to +\infty}J(u_T)+
\limsup_{T \to +\infty}J(v_T)\\
=\lim_{T \to +\infty}\int_{-T}^T J(\chi_{\{u>s\}})\,ds+\lim_{T \to +\infty}\int_{-T}^T J(\chi_{\{v>s\}})\,ds= J(u)+J(v),
\end{multline*}
so that \eqref{sub2} follows.
\par
In order to conclude the proof, we have to check claim \eqref{subadd}.
Let $M:=\max u$. Since $u$ is positive and integer valued, we can write $u=\sum_{i=1}^{M}\chi_{\{u \geq i\}}$.
For any $i \in \{1,\dots,M\}$ we have
$$
\int_{i-1}^i  J(\chi_{\{u>s\}})\,ds= J(\chi_{\{u \geq i\}})
$$
so that
$$
J(u)=\int_0^{+\infty}J(\chi_{\{u>s\}})\,ds=\sum_{i=1}^M
J(\chi_{\{u\geq i\}}).
$$
Hence
$$
J(u)\geq \inf \left\{\sum_{i=1}^m  J(\chi_{E_i}): m \geq 1\,,
u=\sum_{i=1}^m \chi_{E_i}\right\}.
$$
In order to prove the opposite inequality let $u=\sum_{i=1}^m \chi_{E_i}$
for some Borel set $E_i \subseteq \Om$ and $m \geq 1$. Observe that for any $r,s \in \{1,\dots,m\}$ with $r\neq s$ we also have
$$
u=\chi_{E_r \cap E_s}+\chi_{E_r \cup E_s}+\sum_{\begin{subarray}{c} i\neq r\\i\neq
s\end{subarray}}\chi_{E_i}.
$$
From (\ref{sottomodins}) we get
$$
J(\chi_{E_r \cap E_s})+J(\chi_{E_r \cup E_s})+\sum_{\begin{subarray}{c} i\neq r\\i\neq
s\end{subarray}}J(\chi_{E_i}) \leq \sum_{i=1}^{m}  J(\chi_{E_i}).
$$
Then by induction it is easy to see that
\begin{multline*}
\inf\left\{\sum_{i=1}^m J(\chi_{E_i}): m \geq 1\,, u=\sum_{i=1}^m
\chi_{E_i}\right\}
\\
\ge \inf\left\{\sum_{i=1}^m J(\chi_{E_i}): m \geq 1\,, u=\sum_{i=1}^m
\chi_{E_i}\,,E_1 \supseteq E_2 \supseteq \dots \supseteq E_m \right\}=\sum_{i=1}^{M} J(\chi_{\{u\geq i\}})= J(u).
\end{multline*}
Hence claim \eqref{subadd} holds true, so that the proof is concluded.
\end{proof}

The following Proposition deals with the stability of the class $GC(\Om)$ with respect to the $\Gamma$-convergence. 

\begin{proposition}
\label{stgammaconv}
Let $(J_n)_{n \in \N}$ be a sequence of convex functionals in $GC(\Om)$ such that there exists a functional $\tilde J$ defined on characteristic functions of Borel sets
which satisfies the following conditions:
\begin{itemize}
\item[(a)]
for every Borel set $E \subseteq \Om$ and for every sequence of Borel sets $(E_n)_{n \in \N}$ contained in $\Om$ such that
$\chi_{E_n} \to \chi_{E}$ in $L^1(\Om)$ we have
$$
\tilde J(\chi_E)\leq\liminf_{n \to \infty}J_{n}(\chi_{E_n});
$$
\item[(b)] for every Borel set $E \subseteq \Om$ there exists a sequence of Borel sets 
$(E_n)_{n \in \N}$ contained in $\Om$ with $\chi_{E_n} \to \chi_{E}$ in $L^1(\Om)$ and such that
$$
\limsup_{n \to \infty}J_{n}(\chi_{E_n}) \leq
\tilde J(\chi_E).
$$
\end{itemize}
Then setting
$$
J(u):=\int_{-\infty}^{+\infty}\tilde J(\chi_{\{u >s\}})\,ds,
$$
we have $J \in GC(\Om)$ and the sequence $(J_n)_{n \in \N}$ $\Gamma$-converges to $J$ in the $L^1$-topology.
\par
Conversely let $(J_n)_{n \in \N}$ be a sequence of functionals in $GC(\Om)$ which $\Gamma$-converges to a proper functional $J \colon L^1(\Om) \to [0,+\infty]$. Then $J \in GC(\Om)$ and its restriction $\tilde J$ to the family of characteristic functions of Borel subsets of $\Om$ satisfies conditions {\rm(a)} and {\rm(b)}.
\end{proposition}
\begin{remark}\textup{
Notice that it follows that for \textit{convex} functionals in $GC(\Om)$,
the $\Gamma$-convergence is equivalent to the $\Gamma$-convergence
on the corresponding (submodular) set functions,
that is, the restriction to characteristic functions of the original
functionals. However, the last statement in Proposition~\ref{stgammaconv}
is also true without assuming any convexity of the functions $J_n$. Notice that there exist functionals in $GC(\Om)$ which are lower semicontinuous but not convex, so that convexity cannot be gained by relaxation. (It suffices to consider functionals of the form \eqref{Jh} with $\Om$ and $h$ chosen in such a way that the summation involves only one square, and the function $F$ is not submodular on binary vectors.)
}\end{remark}
\begin{proof}[Proof of Proposition~\ref{stgammaconv}]
Notice that $\tilde J$ is, by construction, lower semicontinuous on characteristic functions, so that the map $ s\mapsto \tilde J(\chi_{\{u >s\}})$ is measurable for every $u \in L^1(\Om)$. Hence the definition of $J$ is well posed.
\par
In order to prove the $\Gamma$-convergence result, we need to check $\Gamma$-lim inf and $\Gamma$-lim sup inequalities (see Section \ref{prelsec}). Let us start with the $\Gamma$-lim inf inequality. Let $u_n \to u$ in $ L^1(\Om)$. Up to a subsequence, we can assume that $\chi_{\{u_n>s\}} \to \chi_{\{u>s\}}$ in $L^1(\Om)$ for a.e.~$s \in \R$.
By Fatou's Lemma, the generalized coarea formula \eqref{coarea} and assumption (a) we get
\begin{equation*}
\liminf_{n \to \infty}J_{n}(u_n)\geq
\int_{-\infty}^{+\infty}\liminf_{n \to \infty}J_{n}
(\chi_{\{u_n>s\}})\,ds \geq
\int_{-\infty}^{+\infty}\tilde J(\chi_{\{u>s\}})\,ds=J(u)
\end{equation*}
so that the $\Gamma$-lim inf inequality follows.
\par
Let us come to the $\Gamma$-lim sup inequality. We can clearly assume that the map $s \mapsto \tilde J(\chi_{\{u>s\}})$ belongs to $L^1(\R)$. Notice that the subspace given by (finite) linear combinations of characteristic functions is dense with respect to the energy $J$. In fact, if $0\leq u\leq 1$, by Lemma \ref{appross} we can choose $t \in ]0,1[$ 
such that 
$$
\frac{1}{m}\sum_{\begin{subarray}{c} k\in\mathbb{N}
\end{subarray}}\chi_{\left\{u>\frac{k+t}{m}\right\}} \to u
\qquad\text{ in }L^1(\Om)
$$
and (since $\tilde J(0)=\tilde J(1)=0$)
$$
\lim_{m \to \infty}\frac{1}{m}\sum_{\begin{subarray}{c} k\in\mathbb{N}
\end{subarray}}\tilde J\left(\chi_{\left\{u>\frac{k+t}{m}\right\}}\right)=
\int_{0}^1 \tilde J(\chi_{\{u>s\}})\,ds=
\int_{-\infty}^{+\infty} \tilde J(\chi_{\{u>s\}})\,ds.
$$
The case $m \leq u \leq M$ follows considering the function $(u-m)/(M-m)$. Finally, for $u \in L^1(\Om)$, let us set $u_T:=-T \vee u \wedge T \to u$ for every $T>0$.
Since $u_T \to u$ in $L^1(\Om)$ and $J(u_T) \to J(u)$ as $T \to +\infty$, the density in energy follows by a diagonal argument. By general results of $\Gamma$-convergence, it suffices to prove the $\Gamma$-lim sup inequality for $u$ equal to a linear combination of characteristic functions. Since $J$ is invariant under translation, it is not restrictive to assume that $u=\sum_{k=1}^m a_k \chi_{E_k}$ with $a_k \ge 0$ for every $k=1,\dots,m$, and $E_m \subseteq E_{m-1}\subseteq \dots \subseteq E_1$. In this way we have 
$$
J(u)=\sum_{k=1}^m a_k \tilde J(\chi_{E_k}).
$$
By condition (b), we can find Borel sets $E_k^n$ such that 
$$
\chi_{E_k^n} \to \chi_{E_k}
\qquad\text{in }L^1(\Om)
$$ 
as $n \to \infty$, and
$$
\limsup_{n \to \infty}J_n(\chi_{E_k^n}) \le \tilde J(\chi_{E_k}).
$$
Setting $u_n:=\sum_{k=1}^m a_k\chi_{E_k^n}$, we have $u_n \to u$ in $L^1(\Om)$.
Since $J_n$ is convex and positively one-homogeneous, and hence subadditive, we deduce
\begin{equation*}
\limsup_{n \to \infty}J_n(u_n) \le \sum_{k=1}^m a_k\limsup_{n \to \infty}J_n(\chi_{E_k^n}) \le \sum_{k=1}^m a_k \tilde J(\chi_{E_k})=J(u),
\end{equation*}
so that the $\Gamma$-lim sup inequality is proved.
\par
Finally, the fact that $J \in GC(\Om)$ follows since $J$ and $\tilde J$ coincide on characteristic functions. The proof of the first part of the Proposition is thus complete.
\par
Let us come to the second part. 
Clearly $\tilde J$ satisfies condition {\rm(a)}. In order to prove condition (b), let $E$ be a Borel subset of $\Om$, and let $u_n \in L^1(\Om)$ be such that $u_n \to \chi_E$ in $L^1(\Om)$ and $\limsup_{n \to \infty}J_{n}(u_n)\leq \tilde J(\chi_E)$.
Since for any $\delta \in (0,1)$
$$
J_{n}(u_n) \geq \int_\delta^{1-\delta} J_{n}(\chi_{\{u_n >s\}})\,ds,
$$
there exists $s_n \in (\delta,1-\delta)$ such that
\begin{equation}\label{stima}
J_{n}(\chi_{\{u_n > s_n\}}) \leq \frac{1}{1-2\delta}J_{n}(u_n).
\end{equation}
Let us set $E^\delta_n:=\chi_{\{u_n > s_n\}}$. We have clearly that $\chi_{E^\delta_n} \to \chi_E$ in $L^1(\Om)$ and by \eqref{stima} we deduce
$$
\limsup_{n \to \infty}J_{n}(\chi_{E^\delta_n})\leq \frac{1}{1-2\delta}\limsup_{n \to \infty}J_{n}(u_n)\leq \frac{1}{1-2\delta}\tilde J(\chi_E).
$$
Let us choose now $\delta=1/m$. There exists $n_m$ such that for every $n \ge n_m$ we have 
$$
\|\chi_{E^{1/m}_n}-\chi_E\|_{L^1} \le 1/m
$$
and
$$
J_n(\chi_{E^{1/m}_n}) \le \frac{1}{1-2/m}\tilde J(\chi_E)
$$
Moreover we may assume that $n_m \uparrow \infty$.
If we set $E_n:=E^{1/m}_n$ for $n_m \le n <n_{m+1}$, we have that $(E_n)_{n \in \N}$ is the recovering sequence for which the $\Gamma$-lim sup inequality holds. 
\par
Finally, the fact that $J \in GC(\Om)$ follows now from the first part of the Proposition, and this concludes the proof.
\end{proof}

\section{Discrete approximation of anisotropic total variation}
\label{discrsec}
Let $N\ge 1$ and $\Sigma\subset \Z^N$ be a finite set, and let $F:\{0,1\}^M\to [0,+\infty)$ be a nonnegative submodular function, i.e. $F(u \wedge v)+F(u \vee v) \leq F(u)+F(v)$ for any $u,v \in \{0,1\}^\Sigma$, 
with $F(0)=F(\chi_\Sigma)=0$. 
We extend $F$ to all
vectors $u\in \R^\Sigma$ into a convex function by letting (see Proposition \ref{prop-conv})
\begin{equation}\label{discretecoarea}
F(u)=\int_{-\infty}^{+\infty}F(\chi_{\{u>s\}})\,ds
\end{equation}
where  $\{u>s\}:= \{x\in\Sigma : u(x)>s\}$.
We let
\begin{equation}\label{defrho}
\rho_\Sigma \ :=\ \max_{x\in\Sigma} |x|
\end{equation}
and we assume, in addition, the following coercivity assumption:
\begin{itemize}
\item[{\bf (A)}] $\Sigma$ contains $0$ and
the canonical basis $(e_i)_{i=1}^N$
of $\R^N$, and there exists $c>0$ such that for any
$u \in \R^\Sigma$,
$$
F(u) \geq c\sum_{i=1}^N |u(e_i)-u(0)|.
$$
\end{itemize}
Notice that \eqref{discretecoarea} is a discrete version
of the generalized coarea formula \eqref{coarea}.
\par
Given $h>0$ and $x \in h\Z^N$ let us set
\begin{equation}\label{defQxh}
Q^h_x:=x+h \left[0,1\right[^N.
\end{equation}

Let $V_h$ denote the space of functions $u \colon \R^N \to \R$ such that
$$
u\ =\ \sum_{x\in h\Z^N}u(x)\chi_{Q^h_x}.
$$
Notice that we have $V_h \subseteq L^1_{loc}(\R^N)$.
\par
Let $\Om \subset \R^N$ be an open and bounded set. We denote by $V_h(\Om)$ 
the restriction to $\Om$ of the functions in $V_h$.
Let moreover $I^h(\Om)$ denote the set of $x\in\R^N$
such that $x+h\Sigma\subset\Om$.
We consider the functional
$J_h(\cdot,\Om) \colon L^1(\Om) \to [0,+\infty[$ 
defined as
\begin{equation}
\label{Jh}
J_h(u,\Om):=
\begin{cases}
\displaystyle h^{N-1}\sum_{x \in I^h(\Om)\cap h\Z^N}F(u[x+h\Sigma])
      &\text{if } u \in V_h(\Om)\\
+\infty &\text{if }u \in L^1(\Om)\setminus V_h(\Om),
\end{cases}
\end{equation}
where for any $x\in I^h(\Om)$,
$u[x+h\Sigma]$ is the vector $(u(x+hy)_{y\in \Sigma})$ of $\R^\Sigma$.
\par
The aim of this section is to study the asymptotic behavior of the functionals $J_h(\cdot,\Om)$ as the size mesh $h$ vanishes: it is expected
that they approximate some anisotropic total variation.
 The following Proposition shows
that the functionals $J_h(\cdot,\Om)$ satisfy the generalized coarea
formula \eqref{coarea}.

\begin{proposition}
\label{settdiscr}
The functional $J_h(\cdot,\Om)$ is convex and belongs to $GC(\Om)$.
Moreover, there exist $C_2>C_1>0$ such that
for any open sets $A,B$ with $A\subset\subset
B\subset\subset\Om$, and for any $u\in V_h(\Om)$, we have,
if $h$ is small enough,
\begin{equation}\label{estvariation}
C_1|Du|(A)\ \le\ J_h(u,B)\ \le\ C_2|Du|(\Om)\,.
\end{equation}
\end{proposition}

\begin{proof} From~\eqref{discretecoarea}, we get that
also $J_h$ satisfies~\eqref{coarea}. The submodularity of $F$
yields~\eqref{sottomodins}, hence $J_h$ is convex. 

To show the estimate~\eqref{estvariation}, it is enough
to assume that $u\in V_h(\Om)$ is a characteristic function
(the general case then follows from the coarea formula).
In this case, the left hand side inequality follows
from assumption \textbf{(A)}, while the other follows from
the fact if that $F(u[x+h\Sigma])> 0$ for some $x\in I^h(B)\cap h\Z^N$,
then $u$ takes different values on the set $x+h\Sigma$ so that
its variation on $B(x,\rho_\Sigma h)$ (where $\rho_\Sigma$ is given by~\eqref{defrho})
is at least $h^{N-1}$.
\end{proof}


For every $\nu\in \R^N$ we set
$$
F(\nu\cdot \Sigma):=F\left((\nu\cdot y)_{y\in\Sigma}\right)
$$
The main result of the paper is the following.

\begin{theorem}
\label{main}
Let $\Om\subseteq \R^N$ be a bounded, open and Lipschitz domain, and let
$J_h:=J_h(\cdot,\Om)$ be defined as in \eqref{Jh} for $h>0$.
Then the family $(J_h)_{h>0}$ $\Gamma$-converges in the $L^1$-topology as $h \to 0$
to the functional $J \colon L^1(\Om) \to [0,+\infty]$ given by
\begin{equation}\label{gammalimite}
J(u,\Om)=
\begin{cases}
\displaystyle\int_\Om F\left( \displaystyle\frac{Du}{|Du|}\cdot \Sigma\right)\,d|Du| &\text{if }u \in BV(\Om)\\
+\infty      &\text{if }u \in L^1(\Om)\setminus BV(\Om),
\end{cases}
\end{equation}
where for $u\in BV(\Om)$ the function $Du/|Du|$ stands for the Radon-Nikodym derivative of $Du$ with respect to its total variation $|Du|$.
\end{theorem}

Since $J(\cdot,\Om)$ satisfies the generalized coarea formula (see Proposition \ref{stgammaconv}), we can also write for $u\in BV(\Om)$
\begin{equation}
\label{glimcoarea}
J(u,\Om)=\int_{-\infty}^{+\infty}\Persf(\{u>s\},\Om)\,ds,
\end{equation}
where for any finite-perimeter set $E$ in $\Om$
\begin{equation}
\label{aniper2}
\Persf(E,\Om)\ =\ \int_{\partial^* E\cap \Om}F(\nu_E\cdot \Sigma)\,d\hn.
\end{equation}
In particular 
$$
J(\chi_E,\Om)=\Persf(E,\Om)
$$
for any finite-perimeter set $E$ in $\Om$.
\par
The following Corollary is a consequence of the $\Gamma$-convergence result of the previous Theorem.

\begin{corollary}
\label{cormain}
Let $\Om\subseteq \R^N$ be a bounded open set with Lipschitz boundary, $g \in L^\infty(\Om)$, and let $u_h$ be the solution of 
\begin{equation}
\label{minJn}
\min_{u \in L^1(\Om)} J_h(u,\Om)+\|u-g\|_{L^2(\Om)}^2.
\end{equation}
Then $u_h$ converges in $L^1(\Om)$ for $h \to 0$ to the minimizer $u \in BV(\Om)$ of
\begin{equation}
\label{minJu}
\min_{u \in L^1(\Om)} J(u,\Om)+\|u-g\|_{L^2(\Om)}^2,
\end{equation}
where $J$ is the $\Gamma$-limit of the family $(J_h)_{h>0}$ given by \eqref{gammalimite}.
\end{corollary}

\begin{proof}
Without loss of generality we can suppose $u_h \in V_h(\Om)$ for any $h>0$. Moreover, as a consequence of the submodularity property (iv) in Proposition \ref{prop-1}, we have that the functional $J_h$ decreases by truncation. This entails that $\|u_h\|_\infty \le \|g\|_\infty$ for every $h>0$.
\par
Taking into account~\eqref{estvariation} we get that the total variation of $u_h$ is uniformly bounded. By compactness in $BV$, we deduce that there exists $u \in BV(\Om)$ and a sequence $h_k \to 0$ such that $u_{h_k} \to u$ in $L^1(\Om)$. The convergence is indeed in every $L^p$ for every $1 \le p<+\infty$ since $(u_h)_{h>0}$ is bounded in $L^\infty(\Om)$.
\par
The fact that the limit $u$ is a minimizer of \eqref{minJu}
is a consequence of Proposition~\ref{Gamma-conv-prop}. Since this minimizer
is, in fact, unique, we conclude that
the entire family $(u_h)_{h>0}$ converges to $u$ as $h \to 0$.
\end{proof}

\begin{remark}
\label{Gammalimitemisura}
{\rm
Notice that equality \eqref{glimcoarea} implies that the $\Gamma$-limit $J(u,\Om)$ satisfies \eqref{gammalimite} for every $u \in BV(\Om)$. In fact a direct computation shows that \eqref{gammalimite} holds for simple functions. The extension to the whole $BV(\Om)$ follows by a density argument. Let $u \in BV(\Om)$, and let $(u_n)_{n \in \N}$ be a sequence of simple functions converging in $L^1(\Om)$ to $u$ and such that $|Du_n|(\Om) \to |Du|(\Om)$ as $n \to \infty$. 
From Reshetnyak continuity theorem (see \cite[Thm 2.39]{AFP}) we deduce that
\begin{equation}
\label{convtomeasure}
\lim_{n \to \infty}\int_\Om F\left(\frac{Du_n}{|Du_n|}\cdot\Sigma\right)\,d|Du_n|=
\int_\Om F\left(\frac{Du}{|Du|}\cdot\Sigma\right)\,d|Du|.
\end{equation}
Since $|Du_n|(\Om) \to |Du|(\Om)$, by coarea formula in $BV(\Om)$ we get, up to a subsequence,  $Per(\{u_n>s\},\Om) \to Per(\{u>s\},\Om)$ for a.e. $s \in \R$. Using again Reshetnyak continuity theorem (applied to the measures $\nu\,d\hn \res \partial^*\{u_n>s\}$), we deduce that
\begin{multline*}
\lim_{n \to \infty}\Persf(\{u_n>s\},\Om)=\lim_{n \to \infty}\int_{\partial^*\{u_n>s\}\cap \Om} F(\nu\cdot \Sigma)\,d\hn\\=
\int_{\partial^*\{u>s\} \cap \Om} F(\nu\cdot\Sigma)\,d\hn=
\Persf(\{u>s\},\Om).
\end{multline*}
By the generalized coarea formula for $J$ and the Dominated Convergence Theorem we conclude that
$$
\lim_{n \to \infty}J(u_n,\Om)=J(u,\Om)
$$
so that in view of \eqref{convtomeasure}, the representation \eqref{gammalimite} is proved.
}
\end{remark}

The rest of the Section is devoted to the proof of Theorem \ref{main}. In view of Proposition \ref{stgammaconv} and of Remark \ref{Gammalimitemisura}, in order to study the $\Gamma$-limit of the family $(J_h)_{h>0}$ we can consider the restriction
of $J_h$ to characteristic functions of sets, and show that it
$\Gamma$-converges to the anisotropic perimeter given by~\eqref{aniper2}.

By definition, we need to show that given any sequence $h_m\downarrow 0$,
we have for any Borel set $E\subseteq \Om$:
\begin{itemize}
\item if $\chi_{E_m}\in V_{h_m}(\Om)$ converges to $\chi_E$ in $L^1(\Om)$,
then
\begin{equation}\label{glinf}
\liminf_{m\to\infty} J_{h_m}(\chi_{E_m},\Om)\ \ge\ J(\chi_E,\Om)\,;
\end{equation}
\item there exists a sequence $(E_m)_{m\in \N}$ with $\chi_{E_m}\in V_{h_m}(\Om)$
such that $\chi_{E_m}\to \chi_{E}$ in $L^1(\Om)$ and
\begin{equation}\label{glsup}
\limsup_{m\to\infty} J_{h_m}(\chi_{E_m},\Om)\ \le\ J(\chi_E,\Om)\,.
\end{equation}
\end{itemize}

We prove inequalities \eqref{glinf} and \eqref{glsup} in subsections \ref{secglinf} and \ref{secglsup} respectively. We will use the following ``continuous'' variant of $J_h$, defined on any function and not just on piecewise constant functions of the
class $V_h$: we let, for any $u\in L^1(\Om)$,
\begin{equation}\label{Jch}
J^c_{h}(u,\Om):=\frac{1}{h}\int_{I^h(\Om)}
F(u[x+h\Sigma])\,dx.
\end{equation}
Let $Q_\nu$ is the open unit cube
centered in $0$ with a face orthogonal to $\nu$, and 
\begin{equation}
\label{definu}
I_\nu:=\{x \in \R^N : x \cdot \nu >0\}.
\end{equation}
We have the following result:
\begin{lemma}\label{lemanisotropia}
There holds
\begin{equation}\label{anisotropia}
F(\nu\cdot \Sigma)\ =\ \lim_{h\to 0} J^c_h(\chi_{I_\nu},Q_\nu)\,
\end{equation}
\end{lemma}
\begin{proof}
One has
\[
J^c_h(\chi_{I_\nu},Q_\nu)\ =\ \frac{1}{h}\int_{I^h(Q_\nu)}
F(\chi_{I_\nu}[x+h\Sigma])\,dx.
\]
Now, letting $u(x):=\nu\cdot x$, we have, for any $s\in\R$ and
$x\in \R^N$,
$\chi_{I_\nu}[x+h\Sigma]=\chi_{\{u>s\}}[x+s\nu+h\Sigma]$, so that we may
write
\begin{multline}\label{blah0}
J^c_h(\chi_{I_\nu},Q_\nu)\ =\ \int_{-1/2}^{1/2}
\left(\frac{1}{h}\int_{I^h(Q_\nu)}
F(\chi_{\{u>s\}}[x+s\nu+h\Sigma])\,dx\right)ds
\\ =\ \int_{-1/2}^{1/2}
\left(\frac{1}{h}\int_{I^h(Q_\nu)+s\nu}
F(\chi_{\{u>s\}}[y+h\Sigma])\,dy\right)ds.
\end{multline}
Now, as soon as $|s|< 1/2-\rho_\Sigma h$,
\[
\int_{I^h(Q_\nu)+s\nu}
F(\chi_{\{u>s\}}[y+h\Sigma])\,dy
\ =\ \int_{I^h(Q_\nu)}
F(\chi_{\{u>s\}}[y+h\Sigma])\,dy\,,
\]
and it follows from \eqref{blah0} and the
co-area formula that
\begin{equation}\label{blah1}
J^c_h(\chi_{I_\nu},Q_\nu)\ =\ J^c_h(u,Q_\nu)
\,+\,\epsilon_h
\end{equation}
where the error $\epsilon_h\ge 0$ is
\[
\epsilon_h\ =\ \frac{1}{h}\int_{\frac{1}{2}-\rho_\Sigma h<|s|<\frac{1}{2}}
\left(\int_{I^h(Q_\nu)+s\nu}
F(\chi_{\{u>s\}}[y+h\Sigma])\,dy
- \int_{I^h(Q_\nu)}
F(\chi_{\{u>s\}}[y+h\Sigma])\,dy
\right)\,ds.
\]
One easily checks that $\epsilon_h\le 2(\max_\Sigma F)\rho_\Sigma^2 h\to 0$
as $h\to 0$, and~\eqref{anisotropia} follows from~\eqref{blah1}
and the observation that
\[
J^c_h(u,Q_\nu)\ =\ F(\nu\cdot \Sigma)\frac{|I^h(Q_\nu)|}{|Q_\nu|}
\ \to\ F(\nu\cdot\Sigma)
\]
as $h\to 0$.
\end{proof}

\subsection{$\Gamma$-lim\,inf inequality}
\label{secglinf}
We prove here the inequality \eqref{glinf}. First of all, we may assume (upon extracting
a subsequence) that the lim\,inf in~\eqref{glinf} is a limit, and, also,
that it is finite (otherwise there is nothing to prove), so that in particular
$\sup_{m} J_{h_m}(\chi_{E_m},\Om)=C<+\infty$.
For any $A\subset\subset\Om$, it follows from~\eqref{estvariation}
that for any $m$ large enough,
$ C_1|D\chi_{E_m}|(A)\ \le C$ so that (by lower semicontinuity of the total variation) also $E$ must have finite
perimeter in $\Om$.

We consider the non-negative measures
\[
\mu_m\ =\ h_m^{N-1}\sum_{x\in I^{h_m}(\Om)\cap {h_m}\Z^N}
F(\chi_{E_m}[x+h_m\Sigma]) \delta_x\,,
\]
such that $J_{h_m}(\chi_{E_m},\Om)=\mu_m(\Om)$. Since they are uniformly
bounded, we may also assume that there exists a measure $\mu$
such that $\mu_m\weakst\mu$ as measures. We therefore have
\begin{equation}\label{scimu}
\mu(\Om)\ \le\ \liminf_{m\to\infty} J_{h_m}(\chi_{E_m},\Om)\,,
\end{equation}
hence the thesis follows if we show that
 $\mu\ge F(\nu_E\cdot \Sigma)\hn\res\partial^* E$ as measures.

It is therefore enough to compute the Radon-Nikod\'ym derivative of
the measure $\mu$ with respect to~$\hn\res\partial^* E$, and to
show it is above $F(\nu_E\cdot\Sigma)$.
We know from~\cite[Thm.~5.52]{AFP} that it is given,
for $\hn$-a.e.\ $x\in\partial^* E$, by
\[
\lim_{r\to 0} \frac{\mu(x+rQ_\nu)}{\hn((x+rQ_\nu)\cap \partial^* E)}
\]
for any $\nu\in \R^N$ with unit norm, where $Q_\nu$ is as before the open unit cube
centered in $0$ with a face orthogonal to $\nu$. In particular, at a regular
point $x_0$
of $\partial^\ast E$ we can choose $\nu=\nu_E(x_0)$ (the inner normal
to $E$ at $x_0$) and the limit becomes
\begin{equation}
\ell\ =\ \lim_{r\to 0} \frac{\mu(x_0+rQ_\nu)}{r^{N-1}}.
\end{equation}
Let us now show that $\ell\ge F(\nu\cdot \Sigma)$. We can assume $\ell<+\infty$. Notice that since $x_0$
is regular, we also have
\[
\lim_{r\to 0} \frac{1}{r^N}\int_{x_0+2rQ_\nu} |\chi_{I_\nu}(y-x_0)-\chi_E(y)|\,dy
\ = \ 0.
\]

For a.e.\ $r>0$ (small), we have
\begin{equation*}
\mu(x_0+rQ_\nu)
\ =\ \lim_{m\to\infty} \mu_m(x_0+rQ_\nu)\,,
\end{equation*}
and
\begin{equation*}
\int_{x_0+2rQ_\nu} |\chi_{I_\nu}(y-x_0)-\chi_E(y)|\,dy\ =\ 
\lim_{m\to\infty} \int_{x_0+2rQ_\nu} |\chi_{I_\nu}(y-x_0)-\chi_{E_m}(y)|\,dy\,.
\end{equation*}
Hence, using a diagonal argument, there exist sequences $m_n$ and
$r_n$ such that $h_{m_n}/r_n\to 0$, 
\begin{equation*}
\ell\ =\ \lim_{n\to\infty} \frac{\mu_{m_n}(x_0+r_nQ_\nu)}{r_n^{N-1}}
\end{equation*}
and
\begin{equation*}
\lim_{n\to\infty} \frac{1}{r_n^N}\int_{x_0+2r_n Q_\nu}
|\chi_{I_\nu}(y-x_0)-\chi_{E_{m_n}}(y)|\,dy
\ = \ 0.
\end{equation*}

For each $n$,
we now make the change of variable $y=x_0+r_n z$, and we define
$E'_n= (E_{m_n}-x_0)/r_n \subset (\Om- x_0)/r_n$.
It follows
\begin{equation}\label{ellrisc1}
\ell\ =\ \lim_{n\to\infty} \left(\frac{h_{m_n}}{r_n}\right)^{N-1}
\hspace{-1em}
\sum_{z\in Q_\nu\cap [(h_{m_n}/r_n)\Z^N-x_0/r_n]}
\hspace{-1em}
F(\chi_{E'_n}[z + (h_{m_n}/r_n)\Sigma])
\end{equation}
and
\begin{equation}\label{Enriscalato}
\lim_{n\to\infty} \int_{2Q_\nu}|\chi_{I_\nu}(z)-\chi_{E'_n}(z)|\,dz
\ = \ 0.
\end{equation}
We let $h'_n=h_{m_n}/r_n$ (which goes to $0$), and let
$\theta_n = x_0/h_{m_n}-[x_0/h_{m_n}]$ (the vector whose coordinates
are each the fractional part of the corresponding coordinate of $x_0/h_{m_n}$).
The limit in~\eqref{ellrisc1} becomes:
\begin{equation}\label{ellriscalata}
\ell\ =\ \lim_{n\to\infty} h_n'^{N-1}
\sum_{z\in Q_\nu\cap h'_n(\Z^N-\theta_n)}
F(\chi_{E'_n}[z + h'_n\Sigma]).
\end{equation}
Letting $E''_n:=E'_n+h'_n\theta_n$, we clearly still have 
\begin{equation}\label{Enriscalato2}
\lim_{n\to\infty} \int_{Q_\nu}|\chi_{I_\nu}(z)-\chi_{E''_n}(z)|\,dz
\ = \ 0.
\end{equation}

From now on, to simplify, we will denote $E''_n$ and $h'_n$ by, respectively, $E_n$ and $h_n$. We
consider a basis $(\nu_1,\dots,\nu_N)$ of $\R^N$ such that each $\nu_i$ is orthogonal to a face
of $Q_\nu$, and with $\nu_N=\nu$.
We choose $\eta_0>0$ (small). Writing $x=s\nu_1+x'=(s,x')$
with $x'\cdot \nu_1=0$, we have (using Fubini's Theorem)
\begin{multline}
\int_0^{\eta_0}\left[\frac{1}{h_n}
\int_{1/2-\eta-\rho_\Sigma h_n}^{1/2-\eta}\int_{[-1/2,1/2]^{N-1}}|\chi_{I_\nu}(x)-\chi_{E_n}(x)|\,dx'\,ds\right]d\eta
\\
\ =\ \rho_\Sigma \int_{1/2-\eta_0-\rho_\Sigma h_n}^{1/2}\int_{[-1/2,1/2]^{N-1}}|\chi_{I_\nu}(x)-\chi_{E_n}(x)|\,dx'ds
\end{multline}
which, thanks to \eqref{Enriscalato2}, goes to $0$.

Hence, up to a subsequence, we may assume that for a.e.\ $\eta\in [0,\eta_0]$,
\begin{equation}\label{striptoz}
\frac{1}{h_n}\int_{1/2-\eta-\rho_\Sigma h_n}^{1/2-\eta}
\int_{Q_\nu\cap \nu_1^\perp}|\chi_{I_\nu}(s\nu_1+x')-\chi_{E_n}(s\nu_1+x')|\,dx'ds\ \to\ 0
\end{equation}
as $n\to \infty$. The same holds if we replace in~\eqref{striptoz} $\nu_1$ with any of the $\nu_i$,
$i=2,\dots,N-1$, or if we replace the interval $[1/2-\eta-\rho_\Sigma h_n,1/2-\eta]$ of integration in the
first integral with the interval $[-1/2+\eta,-1/2+\eta+\rho_\Sigma h_n]$. Let us therefore
choose a $\eta>0$ such that all the above mentioned limits are $0$,
and let $Q^\eta_\nu=(1-2\eta)Q_\nu$.

We now extend periodically $E_n\cap Q^\eta_\nu$ in the directions $\nu_1,\dots,\nu_{N-1}$,
into a set $\hat{E}_n\subset \nu^\perp + [-1/2+\eta,1/2-\eta]\nu$,
by letting $\hat{x}\in \hat{E}_n$
if and only if $\hat{x}=x+(1-2\eta)\sum_{i=1}^{N-1} k_i\nu_i$ for some point
$x\in E_n\cap Q^\eta_\nu$ and integers $k_i\in \Z$. It follows
from~\eqref{ellriscalata} that
\begin{equation}\label{ellgejc}
\ell \ \ge\ \liminf_{n\to\infty} J^c_{h_n}(\chi_{\hat{E}_n},Q^\eta_\nu)\,,
\end{equation}
where $J^c_{h_n}$ is as in \eqref{Jch}, and (recalling~\eqref{estvariation}), 
\begin{equation}\label{varbound}
\sup_{n \in \N} |D\chi_{\hat{E}_n}|(\overline{Q^\eta_\nu})
\ \le\ \sup_{n\in \N} |D\chi_{E_n}|(\textrm{int}(Q^\eta_\nu))\,+\,2N
\ <\ +\infty.
\end{equation}

We claim that for any $\tau\in \nu^\perp$, we also have
\begin{equation}\label{jcgejctau}
J^c_{h_n}(\chi_{\hat{E}_n},Q^\eta_\nu) \ \ge\  J^c_{h_n}(\chi_{\hat{E}_n},\tau+Q^\eta_\nu)\,-\,\epsilon_n
\end{equation}
for some error $\epsilon_n\to 0$ which is independent on $\tau$.

Assume to simplify that $\tau=s\nu_1$ for some $s\in\R$. If $s$ is $(1-2\eta)$ times
an integer, then~\eqref{jcgejctau}
is obvious. If not, we may assume without loss of generality that $0<s<1-2\eta$. We have
\[
 J^c_{h_n}(\chi_{\hat{E}_n},s\nu_1+Q^\eta_\nu)\ =\ 
\frac{1}{h_n}\int_{I^{h_n}(s\nu_1+Q^\eta_\nu)} F(\chi_{\hat{E}_n}[x+h_n\Sigma])\,dx.
\]
The domain of integration is split into three parts $D_1=I^{h_n}(s\nu_1+Q^\eta_\nu)\cap I^{h_n}(Q^\eta_\nu)$,
$D_2=I^{h_n}(s\nu_1+Q^\eta_\nu)\cap I^{h_n}((1-2\eta)\nu_1+Q^\eta_\nu)$, and
$D_3= I^{h_n}(s\nu_1+Q^\eta_\nu)\setminus (D_1\cup D_2)$. Since
$D_1$ and $(D_2-(1-2\eta)\nu_1)$ are disjoint subsets of $I^{h_n}(Q^\eta_\nu)$, it follows
\begin{equation}\label{boundhatuno}
 J^c_{h_n}(\chi_{\hat{E}_n},s\nu_1+Q^\eta_\nu)\ \le\ J^c_{h_n}(\chi_{\hat{E}_n},Q^\eta_\nu)
\,+\, \frac{1}{h_n}\int_{D_3} F(\chi_{\hat{E}_n}[x+h_n\Sigma])\,dx.
\end{equation}
We have $D_3\subset S_1^+\cup ((1-2\eta)\nu_1+S_1^-)$, where for $i=1,\dots,N-1$,
\[
S_i^\pm \ =\ 
\left\{ x =\sum_{j=1}^N x_j\nu_j\,:\, \frac{1}{2}-\eta - \rho_\Sigma h_n\le \pm x_i < \frac{1}{2}-\eta \,, 
\ |x_j| < \frac{1}{2}-\eta\ \forall j\neq i\right\}.
\]

Let us show that $(1/h_n)\int_{S_1^+}  F(\chi_{\hat{E}_n}[x+h_n\Sigma])\,dx\to 0$ as $n\to\infty$. We have,
using the notation $x=s\nu_1+x'$, $x'\cdot \nu_1=0$, and the change of variable $s=1/2-\eta-h_n\xi$,
\begin{equation}\label{cambiovariabilestrip}
\frac{1}{h_n}\int_{S_1^+}  F(\chi_{\hat{E}_n}[x+h_n\Sigma])\,dx\ =\ \int_0^{\rho_\Sigma}
\int_{[-1/2+\eta,1/2-\eta]^{N-1}}F(v_n(\xi,x'))\,dx'\,d\xi
\end{equation}
where for each $\xi,x'$, $v_n(\xi,x')\in \{0,1\}^\Sigma$  is the vector
$\chi_{\hat{E}_n}[(1/2-\eta-h_n\xi)\nu_1+x'+h_n\Sigma]$. We observe that from~\eqref{striptoz}, we have
(using the same change of variable, and observing that $\chi_{I_\nu}$ depends only on $x'_N$)
\begin{equation*}
\int_{0}^{\rho_\Sigma }\int_{[-1/2,1/2]^{N-1}}
|\chi_{I_\nu}(x'_N)-\chi_{E_n}(1/2-\eta-h_n\xi,x')|\,dx'd\xi \ \to\ 0
\end{equation*}
as $n\to\infty$. In particular, up to a subsequence, we may assume that each component of the
vector $v_n(\xi,x')$ converges to $\chi_{I_\nu}(x'_N)$ as $n\to\infty$, for a.e.\ $(\xi,x')$.
Notice that these components take only the value $0$ or $1$, hence, they must be equal to
$\chi_{I_\nu}(x'_N)$ for $n$ large enough.
Since $F(0)=F(\chi_\Sigma)=0$, it follows that $F(v_n(\xi,x'))\to 0$ a.e., and since $F$ is
bounded we find (using Lebesgue's theorem) that the integral in the right-hand side
of~\eqref{cambiovariabilestrip}
goes to $0$.

In the same way, we can show that
$(1/h_n)\int_{S_1^-}  F(\chi_{\hat{E}_n}[x+h_n\Sigma])\,dx\to 0$, so that~\eqref{jcgejctau} follows
from~\eqref{boundhatuno} and the inclusion $D_3\subset S_1^+\cup ((1-2\eta)\nu_1+S_1^-)$. 
In the general case (if $\tau$ is not parallel to $\nu_1$), we can
show in the same way that~\eqref{jcgejctau} holds provided
\[
\epsilon_n\ =\ \frac{1}{h_n} \sum_{i=1}^{N-1}
\int_{S_i^-\cup S_i^+} F(\chi_{\hat{E}_n}[x+h_n\Sigma])\,dx\,;
\]
our choice of $\eta$ guarantees that $\epsilon_n$ still goes to zero.

Let $u_n \colon Q^\eta_\nu \to [0,1]$ be the average of $\chi_{\hat{E}_n}$ on
each hyperplane orthogonal to $\nu$, given by
$$
u_n(x):=\frac{1}{(1-2\eta)^{N-1}}\int_{Q^\eta_\nu\cap\nu^\perp}
\chi_{\hat{E}_n}(x+x')\,dx'\,.
$$
It is clear that $u_n$ depends only on $x\cdot \nu$.
Since
$$
\int_{Q^\eta_\nu}|u_n-\chi_{I_\nu}|\,dx \le \int_{Q^\eta_\nu}|\chi_{\hat{E}_n}-\chi_{I_\nu}|\,dx,
$$ 
we deduce
$$
u_n \to \chi_{I_\nu} \quad \text{in }L^1(Q^\eta_\nu).
$$
Notice that for a.e.~$x\in Q^\eta_\nu$, the vector $u_n[x+h_n\Sigma]$ is also the average
on the hyperplane through $x$ orthogonal to $\nu$ of the vectors $\chi_{\hat{E}_n}[\cdot+h_n\Sigma]$,
so that, by Jensen's inequality and using~\eqref{jcgejctau}
\begin{equation*}
J^c_{h_n}(\chi_{\hat{E}_n},Q^\eta_\nu) \ \geq\ J^c_{h_n}(u_n,Q^\eta_\nu) \ -\epsilon_n .
\end{equation*}

Together with~\eqref{ellgejc}, it yields
\begin{equation}\label{ellgehomog}
\ell \ \ge\ \liminf_{n\to\infty} J^c_{h_n}(u_n,Q^\eta_\nu)\,.
\end{equation}

It is clear that~\eqref{varbound} also yields a uniform bound on the total variations
$|D u_n|(Q^\eta_\nu)$, $n\in \N$.
Let us fix $\varepsilon, \delta \in ]0,1[$.
By the generalized coarea formula for $J^c_{h_n}$, and the coarea formula in $BV$ we get
$$
\begin{aligned}
C \geq J^c_{h_n}(u_n,Q^\eta_\nu)+\varepsilon |D u_n|(Q^\eta_\nu) & =\int_0^1 [J^c_{h_n}(\chi_{\{u_n>s\}},Q^\eta_\nu)+\varepsilon \hn(\partial^*\{ u_n>s\}\cap Q^\eta_\nu)]\,ds\\
& \geq \int_{\delta}^{1-\delta} [J^c_{h_n}(\chi_{\{u_n>s\}},Q^\eta_\nu)+\varepsilon \hn(\partial^*\{u_n>s\}\cap Q^\eta_\nu)]\,ds
\end{aligned}
$$
for some positive constant $C$. We deduce that there exists $s_n \in ]\delta,1-\delta[$ such that
\begin{equation}\label{maggJhk}
J^c_{h_n}(u_n,Q^\eta_\nu)+\varepsilon |D u_n|(Q^\eta_\nu)
\geq (1-2\delta)
J^c_{h_n}(\chi_{\{u_n>s_n\}},Q^\eta_\nu)
\end{equation}
and
\begin{equation}\label{maggper}
C \geq \varepsilon (1-2\delta) \hn(\partial^*\{u_n>s_n\} \cap Q^\eta_\nu).
\end{equation}
We have
$$
v_n:=\chi_{\{u_n>s_n\}} \to \chi_{I_\nu}
\qquad\text{in }L^1(Q^\eta_\nu).
$$ 
By \eqref{maggper} and since $u_n$ depends only on $x\cdot \nu$,
we deduce that there exists, up to a subsequence,
an odd number $M$ independent of $n$ such that
$$
\partial^* \{v_n=1\}\ =\ Q^\eta_\nu\,\cap \,\bigcup_{k=1}^M (a_n^k\nu+ \nu^\perp)
$$
with 
$$
-\frac{1}{2}+\eta \,<\,a_n^1\,<\,a_n^2\,<\,\dots\,<\,a_n^M \,<\, \frac{1}{2}-\eta.
$$
Without loss of generality, we can assume that $a_n^k \to 0$ for every $k=1,\dots,M$, 
$v_n=0$ if $x \cdot \nu < a_n^1$ and $v_n=1$ if $x \cdot \nu >a_n^M$. Indeed, if
some of the points $a_n^k$ do not go to zero as $n\to\infty$, we can lower the
energy $J^c_{h_n}(v_n,Q^\eta_\nu)$ by ``removing'' from $v_n$ the corresponding discontinuities.

\par
If $M>1$, let us consider the function $\tilde v_n = \chi_{\{x\cdot\nu > a_n^2\}}$.
Since $v_n\vee\tilde{v}_n$ is both a translate of $\tilde{v}_n$ (of $(a_n^1-a_n^2)\nu$)
and of $\chi_{I_\nu}$ (of $a_n^2\nu$),
for $n$ large we have
$$
J^c_{h_n}(\tilde v_n,Q^\eta_\nu)
\ =\ J^c_{h_n}(v_n \vee \tilde v_n,Q^\eta_\nu)\ =\ J^c_{h_n}(\chi_{I_\nu},Q^\eta_\nu).
$$
By the submodularity property of $F$, we obtain setting $v^1_n:=v_n \wedge \tilde v_n$
that for $n$ large enough,
$$
J^c_{h_n}(v^1_n,Q^\eta_\nu) \ \leq\ J^c_{h_n}(v_n,Q^\eta_\nu).
$$
If $M=3$ then $v^1_n=\chi_{\{x\cdot\nu>a^3_n\}}$ is a translate of $\chi_{I_\nu}$.
If $M>3$, then we can reiterate this argument, replacing now $a_n^2$ with $a_n^4$:
after a finite number of steps we eventually find a
translate of $\chi_{I_\nu}$ with energy lower than $J^c_{h_n}(v_n,Q^\eta_\nu)$.
Hence, for $n$ large enough,
\begin{equation}
\label{maggfinale}
J^c_{h_n}(\chi_{I_\nu},Q^\eta_\nu)\ \leq\ J^c_{h_n}(v_n,Q^\eta_\nu).
\end{equation}

By \eqref{anisotropia} and a straightforward scaling argument, the left-hand side
of~\eqref{maggfinale} goes to $(1-2\eta)^{N-1}F(\nu\cdot\Sigma)$ as $n\to\infty$. Recalling~\eqref{ellgehomog}
and~\eqref{maggJhk}, we deduce
\[
\ell+\e\sup_{n\in \N}|Du_n|(Q^\eta_\nu)\ \ge\ (1-2\delta)(1-2\eta)^{N-1}F(\nu\cdot\Sigma)\,,
\]
and since $\e$, $\delta$ and $\eta$ can be chosen arbitrarily small we deduce that
$\ell\ge F(\nu\cdot\Sigma)$. It follows that $\mu\ge F(\nu\cdot\Sigma)\hn\res \partial^* E$,
which together with~\eqref{scimu} yields~\eqref{glinf}.

\subsection{$\Gamma$-lim\,sup inequality}\label{secglsup}
Let us now show the inequality \eqref{glsup}. We first show the following generalization of the formula \eqref{anisotropia}:
\begin{lemma}\label{lemanisotropiaA}
Let $A\subset \Om$ be an open set, and assume $\hn(\partial I_\nu\cap \partial A)=0$.
Then
\begin{equation}\label{anisotropiaA}
\lim_{h\to 0} J^c_h(\chi_{I_\nu},A)\ =\ \hn(\partial I_\nu\cap A)F(\nu\cdot\Sigma)\,.
\end{equation}
\end{lemma}

\begin{proof}
Let $\e>0$. Let $(\nu_i)_{i=1}^N$ be a basis of $\R^N$, with $\nu_N=\nu$,
and let us consider the family $\mathcal{Q}_\e$ of all cubes $Q\subset A$,
of side $\e$,
centered at the points $\sum_{i=1}^N \e k_i\nu_i$ for $k_i\in \Z^N$,
and such that each face is orthogonal to some $\nu_i$.

We have obviously
\[
\lim_{h\to 0} J^c_h(\chi_{I_\nu},A)\ \ge\ \sum_{Q\in\mathcal{Q}_\e}
 \liminf_{h\to 0} J^c_h(\chi_{I_\nu},Q)
\]
while (by~\eqref{anisotropia} and a simple scaling argument) for each $Q$,
\[
\lim_{h\to 0} J^c_h(\chi_{I_\nu},Q)\ =\ \e^{N-1}F(\nu\cdot\Sigma)\,.
\]
Hence,
\[
\lim_{h\to 0} J^c_h(\chi_{I_\nu},A)\ \ge\ \hn\left(\partial I_\nu\cap
\bigcup_{Q\in\mathcal{Q}_\e}Q\right)F(\nu\cdot\Sigma)\,,
\]
letting then $\e\to 0$, it shows ``$\ge$'' in~\eqref{anisotropiaA}.

To show the reverse inequality, we let, for each $\e>0$,
$A^\e=\bigcup_{Q\in\mathcal{Q}_\e} Q$, and first observe that
for any $\eta>0$,
\[
J^c_h(\chi_{I_\nu},A^\e)\ \le\ \sum_{Q\in\mathcal{Q}_\e}
J^c_h(\chi_{I_\nu},(1+\eta)Q)
\]
as soon as $h<\eta\e/(2\rho_\Sigma)$, where here $(1+\eta)Q$ denotes the
cube of same center as $Q$, and dilated by the factor $1+\eta$.
Taking the limit $h\to 0$, we find
\begin{equation}\label{limsupA}
\limsup_{h\to 0} J^c_h(\chi_{I_\nu},A^\e)\ \le\ (1+\eta)^{N-1}
\hn(\partial I_\nu\cap A^\e)F(\nu\cdot\Sigma)
\ \le\ \hn(\partial I_\nu\cap A)F(\nu\cdot\Sigma)\,+\,C\eta\,.
\end{equation}
for any $\eta>0$. Now,
\[
J^c_h(\chi_{I_\nu},A)\ =\ J^c_h(\chi_{I_\nu},A^\e)\,+
\,\frac{1}{h}\int_{I^h(A)\setminus I^h(A^\e)}F(\chi_{I_\nu}[x+h\Sigma])\,dx
\]
so that~\eqref{anisotropiaA} will follow from~\eqref{limsupA} if we show that
\begin{equation}\label{errorstrip}
\frac{1}{h}\int_{I^h(A)\setminus I^h(A^\e)}F(\chi_{I_\nu}[x+h\Sigma])\,dx\ \to\ 0
\end{equation}
as $h\to 0$ and $\e \to 0$. The integrand above is zero when $x$ is at distance larger
than $h\rho_\Sigma$ to the interface $\partial I_\nu$, and $x\in A$ is
out of the domain of integration as soon as it is at distance larger than $h\rho_\Sigma$
to $A\setminus A^\e$, for instance
when $x\in A_{2\e}:=\{\xi \in A\,:\,{\rm dist}(\xi,\partial A)>2\e\}$. Hence, the error~\eqref{errorstrip} is bounded by
\begin{equation}\label{boundsenzaG}
(\max_\Sigma F)\frac{|(A\setminus A_{2\e})\cap (\partial I_\nu+[-h\rho_\Sigma,
h\rho_\Sigma]\nu)|}{h}\,.
\end{equation}
Let $G\subset \partial I_\nu$ be a relatively open set which contains
$\partial I_\nu\cap\partial A$. We claim that there exists $\delta_0>0$ small
such that the projection of $(A\setminus A_\delta)\cap (\partial I_\nu+[-\delta,\delta]\nu)$
onto $\partial I_\nu$ is contained in $G$ for any $\delta<\delta_0$: 
if not, one finds a sequence $\delta_n\to 0$ and points
$x_n \in (A\setminus A_{\delta_n})\cap (\partial I_\nu+[-\delta_n,\delta_n]\nu)$ which
project outside of $G$, but then, any limit point of this sequence should
be in $\partial A\cap\partial I_\nu$ but outside $G$, which is not possible.
\par
Then, if we choose $\e<\delta_0/2\sqrt{N}$, and $h$ small enough, we have
\[
(A\setminus A_{2\e})\cap (\partial I_\nu+[-h\rho_\Sigma,
h\rho_\Sigma]\nu)\ \subset\ G+[-h\rho_\Sigma,h\rho_\Sigma]\nu
\]
so that the lim\,sup of~\eqref{boundsenzaG}, as $h\to 0$, is less than
\[
\rho_\Sigma \left(\max_\Sigma F\right)\,\hn(G).
\]
Since we assumed that $\hn(\partial I_\nu\cap\partial A)=0$, $\hn(G)$ may
be chosen as small as we wish so that~\eqref{errorstrip} holds. Hence the
Lemma is proved. 
\end{proof}

Now, we show the following estimate:
\begin{lemma}\label{lemminksup}
Let $A\subset \Om$ be an open set with continuous boundary, and let $E\subset \Om$ be a finite-perimeter set such that $\hn(\overline{\partial^* E\cap A}\setminus\partial^* E\cap A)=0$.
Then,
\begin{equation}\label{minksup}
\limsup_{h\to 0} J^c_h(\chi_E,A)\ \le\ (2\rho_\Sigma \max_\Sigma F)\,
\hn(\partial^* E\cap A)\,.
\end{equation}
\end{lemma}

\begin{proof}
Since $A$ has a continuous boundary, for $h$ small enough, $x+h\Sigma$ contains
points of Lebesgue density of $E$ both zero and one
only if $x$ is close enough to $\partial^* E\cap A$, namely,
$\dist(x,\partial^* E\cap A)\le h\rho_\Sigma$.
Hence,
\[
J^c_h(\chi_E,A)\ =\ \frac{1}{h} \int_{I^h(A)} F(\chi_E[x+h\Sigma])\,dx\ \le
\ \frac{|\{x\in A\,:\,\dist(x,\partial^* E\cap A)\le h\rho_\Sigma\}|}
{h}(\max_\Sigma F).
\]
By standard results on the Minkowski contents~\cite[Thm 3.2.39]{Federer}, the last fraction
goes to \linebreak 
$2\rho_\Sigma\hn(\partial^* E\cap A)$ as $h\to 0$, which shows the
Lemma.
\end{proof}

We deduce the following:

\begin{corollary}\label{limpoly}
Let $\Om\subset\R^N$ have a continuous boundary, and let $E\subset \Om$ be a set whose boundary is made of a finite union of subsets of $(N-1)$-dimensional hyperplanes. Then
\[
\lim_{h\to 0} J^c_h(\chi_E,\Om)\ =\ \int_{\partial E\cap\Om}F(\nu_E(x)\cdot\Sigma)\,d\hn(x),
\]
where $\nu_E(x)$ is the inner normal to $E$ at $x$.
\end{corollary}

\begin{proof}
By assumption, $\partial E\cap\Om=\bigcup_{i=1}^M P_i$ where
$P_i\subset (x_i+\partial I_{\nu_i})$ for some $x_i\in \R^N$
and $\nu_i\in \R^N$ with unit norm, moreover we assume that $\nu_i=\nu_E$
on $P_i$ ($\nu_i$ points towards the interior of $E$).
(We also assume that the $P_i$ are ``maximal'', in the sense that
$\overline{P}_i\cap\overline{P}_j\neq \emptyset\Rightarrow \nu_i\neq \nu_j$ for
any $i\neq j$.)
 Let $S=\Om\cap\bigcup_{i=1}^M
\partial P_i$, where the $\partial P_i$ is the relative boundary of the face
inside the hyperplane $(x_i+\partial I_{\nu_i})$: it is a $(N-2)$-dimensional
set, with $\mathcal{H}^{N-2}(S)<+\infty$.

We choose $\e>0$ and build a finite covering $(A_i)_{i=1}^{M+2}$
of $\partial E\cap\Om$ with bounded open sets with continuous boundary, as follows:
$A_{M+1}=\{x\in\Om\,:\, \dist(x,\partial \Om)<\e\}$,
$A_{M+2}=\{x\in \Om,:\, \dist(x,S)<\e\}$, and
$A_i$ is a neighborhood of $P_i\setminus (A_{M+1}\cup A_{M+2})$ which
does not intersect $P_j$ for $j\neq i$. By Lemmas~\ref{lemanisotropiaA} and~\ref{lemminksup},
\[
\left|\lim_{h\to 0} J^c_h(\chi_E,\Om)\,-\,\sum_{i=1}^M \hn(P_i\cap A_i)F(\nu_i\cdot\Sigma)\right|
\,\le\,
C\left(\hn(\partial E\cap A_{M+1})+\hn(\partial E\cap A_{M+2})\right)
\]
with $C=2\rho_\Sigma\max_\Sigma F$. Letting $\e\to 0$ shows the Corollary.
\end{proof}

We are now able to show the following Proposition, which essentially
shows the $\Gamma$-convergence of $J_h$ to $J$ on the restricted class
of polyhedral sets.

\begin{proposition}\label{limsupcrystal}
Let $\Om\subset\R^N$ have a continuous boundary, and let $E\subset \Om$ be a set whose boundary is made of a finite union
of subsets of $(N-1)$-dimensional hyperplanes. Then, there exist sets
 $E_h$ with $\chi_{E_h}\in V_h(\Om)$, $\chi_{E_h}\to \chi_E$ in $L^1(\Om)$
as $h\to 0$, and
\[
\lim_{h\to 0} J_h(\chi_{E_h},\Om)\ =\ \int_{\partial E\cap\Om}F(\nu_E(x)\cdot\Sigma)\,d\hn(x)
\ =\ J(\chi_E,\Om)\,.
\]
\end{proposition}

\begin{proof}
We have, making in~\eqref{Jch} the change of variable $x=z+y$ with $z\in h\Z^N$
and $y\in [0,h[^N$,
\[
J^c_h(\chi_E,\Om)\ =\ \frac{1}{h^N}\int_{[0,h[^N}
\left[h^{N-1}\hspace{-5mm}
\sum_{z\in (I^h(\Om)-y)\cap h\Z^N} F(\chi_E[y+z+h\Sigma])\right]dy\,
\]
so that (by Corollary~\ref{limpoly}) we can choose for each $h$ a $y_h\in [0,h[^N$ such that
\begin{equation}\label{supyh}
\limsup_{h\to 0} \ h^{N-1}\hspace{-5mm}
\sum_{z\in (I^h(\Om)-y_h)\cap h\Z^N} F(\chi_E[y_h+z+h\Sigma])
\ \le\  J(\chi_E,\Om)\,,
\end{equation}
we can assume moreover that no point in $y_h+h\Z^N\cap \Om$ lies on $\partial E$.

For each $\e>0$ we let $\Om_\e=\{x\in \Om\,:\,\dist(x,\partial \Om)>\e\}$.
We fix $\e>0$, and define $E^\e_h$ as follows (for $h$ small enough):
\[
\chi_{E^\e_h}\ =\ \sum_{z\in \Om_\e\cap h\Z^N} \chi_{E}(y_h+z)\chi_{Q_z^h}
\,+\,\sum_{z\in (\Om\setminus \Om_\e)\cap h\Z^N} \chi_{E}(z)\chi_{Q_z^h}
\]
where $Q_z^h$ is the cube defined in~\eqref{defQxh}. It is not difficult
to show that $\chi_{E^\e_h}\to \chi_E$ in $L^1$ as $h\to 0$, in fact,
it converges locally uniformly in $\Om\setminus\partial E$.
We have if $h$ is small enough
\begin{multline}
\label{splitJh}
J_h(\chi_{E^\e_h},\Om)\ \le\ h^{N-1}\sum_{x \in I^h(\Om_\e) \cap h\Z^N} F(\chi_{E^\e_h}[x+h\Sigma])+
h^{N-1}\sum_{x\in I^h(\Om\setminus \Om_{2\e})\cap h\Z^N} F(\chi_{E^\e_h}[x+h\Sigma])
\\
\le
h^{N-1}\hspace{-5mm}
\sum_{z\in (I^h(\Om)-y_h)\cap h\Z^N} F(\chi_E[y_h+z+h\Sigma])+
h^{N-1}\sum_{x\in I^h(\Om\setminus \Om_{2\e})\cap h\Z^N} F(\chi_{E^\e_h}[x+h\Sigma])
\end{multline}
Now, thanks to~\eqref{supyh} we get
\begin{equation}
\label{splitJhbis}
\limsup_{h \to 0}J_h(\chi_{E^\e_h},\Om) \le J(\chi_E,\Om)
+\limsup_{h \to 0}h^{N-1}\hskip-5mm\sum_{x\in I^h(\Om\setminus \Om_{2\e})\cap h\Z^N} F(\chi_{E^\e_h}[x+h\Sigma]).
\end{equation}
In the sum, on the other hand, $F(\chi_{E^\e_h}[x+h\Sigma])$ is not zero only
when some point of $x+h\Sigma$ lies in $E^\e_h$ and some other in $\Om\setminus E^\e_h$,
and such $x$ are at distance at most $h(1+\rho_\Sigma)$ from $\partial E\cap \Om$,
so that
\[
Q_x^h\ \subset\ \{ \xi\in \Om\setminus\Om_{2\e+h\sqrt{N}}\,
:\,\dist(\xi,\partial E\cap\Om)\le h(1+\sqrt{N}+\rho_\Sigma)\}
\]
Therefore, the sum is bounded by
\[
\frac{|\{ \xi\in \Om\setminus \Om_{2\e+h\sqrt{N}}\,
:\,\dist(\xi,\partial E\cap\Om)\le h(1+\sqrt{N}+\rho_\Sigma)\}|}{h}\,,
\]
which goes to $2(1+\sqrt{N}+\rho_\Sigma)\hn(\partial E\cap(\Om\setminus\Om_{2\e}))$
as $h\to 0$. We deduce from~\eqref{splitJhbis} that
\[
\limsup_{h\to 0} J_h(\chi_{E^\e_h},\Om)\ \le\ J(\chi_E,\Om)\,+\,
C\hn(\partial E\cap(\Om\setminus\Om_{2\e}))\,.
\]
Since $\hn(\partial E\cap(\Om\setminus\Om_{2\e}))\to 0$
as $\e\to 0$, using a diagonal argument, we deduce the thesis of the Corollary.
\end{proof}

\begin{corollary}\label{glimcrystal}
Let $\Om\subset\R^N$ have a continuous boundary, and let $E$ be a polyhedral set
in $\Om$ in the previous sense. Then,
\[
(\Gamma-\lim_{h\to 0} J_h(\cdot,\Om))(\chi_E)\ =\ J(\chi_E,\Om)\,.
\]
\end{corollary}
\begin{proof}
It follows from~\eqref{glinf} (which has been shown in Subsection~\ref{secglinf}) and
from Proposition~\ref{limsupcrystal}. 
\end{proof}

Now, we are in a position to show that~\eqref{glsup} holds.

\begin{proposition}\label{above}
Let $\Om$ be a bounded open set of $\R^N$ with Lipschitz boundary, 
and let $E\subset\Om$ be a set with finite perimeter in $\Om$.
Then for every $h>0$ there exists $E_h$ with $\chi_{E_h} \in V_h(\Om)$,
such that $\chi_{E_h}\to\chi_{E}$ in $L^1(\Om)$ as
$h\to 0$ and
\begin{equation}
\label{mainestabove}
\limsup_{h \to 0}J_h(\chi_{E_h},\Om)\ \leq\ J(\chi_{E},\Om).
\end{equation}
\end{proposition}

\begin{proof}
Since $\Om$ is Lipschitz, we can extend $E$ outside of $\Om$ into a finite-perimeter
subset of $\R^N$ (still denoted $E$) such that $|D\chi_E|(\partial\Om)=0$.
Then, standard approximation 
arguments show that there exists a sequence of polyhedral sets $G^n$ such that
$\chi_{G^n}$ converges to $\chi_{E}$ strongly in $L^1(\R^N)$, and with
\[
\lim_{n\to\infty} |D\chi_{G^n}(\R^N)|\ =\ |D\chi_{E}|(\R^N)\,.
\]
This can be seen, for instance, by approximating $\chi_{E}$ by smooth
functions (by convolution) and then approximating these smooth functions
by piecewise linear functions, such as ``P1'' finite-elements. Then, an
appropriate thresholding of these functions provides the sequence $G^n$.
The Reshetnyak continuity Theorem (see section \ref{prelsec}), and
$|D\chi_E|(\partial\Om)=0$ yield
\[
\lim_{n\to\infty} J(\chi_{G^n},\Om)\ =\ 
J(\chi_E,\Om)\,.
\]
By Proposition~\ref{limsupcrystal}, we know that for each $n$ there exists $G^n_h$
converging to $G^n$, such that $\limsup_{h \to 0} J_h(G^n_h,\Om)\le J(\chi_{G^n},\Om)$.
We construct the family $(E_h)_{h>0}$ from the $G_h^n$, by a diagonal argument as follows. For every $n$ there exists $h_n$ such that $h_n \downarrow 0$ as $n \to \infty$ such that for every $h \le h_n$ we have $\|\chi_{G^n_h}-\chi_{G^n}\|_{L^1(\Om)} \le 1/n$ and $J_h(\chi_{G^n_h},\Om)\le J(\chi_{G^n},\Om)+1/n$. If we set $E_h:=G^n_h$ for $h_{n+1}<h\le h_n$, the result follows.
\end{proof}

\section{Examples}

Let us describe a few cases. First of all, the standard nearest-neighbor
interaction on a square grid corresponds to the situation where
$\Sigma = \{(0,0),(1,0),(0,1)\}$ and, for $u\in \R^\Sigma$,
$F(u) = |u(1,0)-u(0,0)|+|u(0,1)-u(0,0)|$. It is obvious, in this
case (as in any other case where $F$ is a sum of pair interactions)
that the $\Gamma$-limit of $J_h$ is the anisotropic total
variation given by~\eqref{totvarintr}, in this case, $\int_\Om |Du|_1$ where 
$|p|_1=|p_1|+|p_2|$ is the $1$-norm in $\R^2$.

\begin{figure}[thb!]
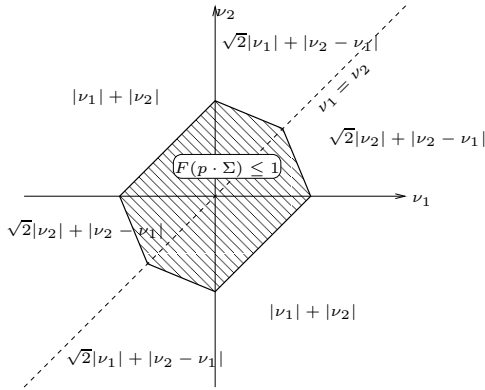

\include{example1} 
\caption{Values of $F$ given by~\eqref{eqexample1}}
\label{figexample1}
\end{figure}

Less trivial situations are when $F$ cannot be reduced to a sum
of pair interactions, such as, still with the same set $\Sigma$,
the functions $F$ defined by $F(0_\Sigma)=F(1_\Sigma)=0$,
$F(1_\Sigma-u)=F(u)$ for any $u\in \{0,1\}^\Sigma$, and
\[
F\left(\begin{array}{ll} 0 \\ 0 & 1\end{array} \right)\,=\,1\,,
\qquad F\left(\begin{array}{ll} 1 \\ 0 & 0\end{array}\right)\,=\, 1\,,
\qquad F\left(\begin{array}{ll} 1 \\ 0 & 1\end{array}\right)\,=\,\sqrt{2}\,.
\]
This $F$ can also be checked to be submodular. Now,
the limit density is given by
\begin{equation}\label{eqexample1}
\nu=(\nu_1,\nu_2) \ \mapsto\ 
F\left(\begin{array}{ll} \nu_2 \\ 0 & \nu_1\end{array} \right)\,,
\end{equation}
see Figure~\ref{figexample1} for the expression of $F$, and where
we also have plotted the shape of the ``Wulff shape'' $F(p\cdot\Sigma)\le 1$.
Notice that in this case, we have 
chosen $F(\theta) = \sqrt{(\theta_{1,0}-\theta_{0,0})^2+
(\theta_{0,1}-\theta_{0,0})^2}$ for $\theta \in \Sigma^{\{0,1\}}$:
however, clearly, we get a limit energy which is not
$\int_\Om |Du|$ (which would be the $\Gamma$-limit of
$u\mapsto h\sum_{i,j} \sqrt{(u_{i+1,j}-u_{i,j})^2+
(u_{i,j+1}-u_{i,j})^2}$) but an anisotropic (crystalline) total variation.
This is of course due to the fact that the latter does not satisfy
the generalized coarea formula.

\section*{Acknowledgments}
A.\,C. is supported by ANR, program ``MICA'', grant ANR-08-BLAN-0082
and by the CNRS. A.\,G. is supported by the Italian Ministry of University and Research, project ``Variational problems with multiple scales'' 2006.  
L.\,L. was partially supported by
the Istituto Nazionale di Alta Matematica (Roma, Italy) and by the Universit\`{a} Cattolica
del Sacro Cuore (Brescia, Italy) during the permanence of the author at the Centre de Math\'{e}matiques Appliqu\'{e}es of the \'{E}cole Polytechnique (Palaiseau, France).

\bibliographystyle{plain}
\bibliography{CGL}











\end{document}

%% file: example1.tex
\begin{picture}(0,0)%
\includegraphics{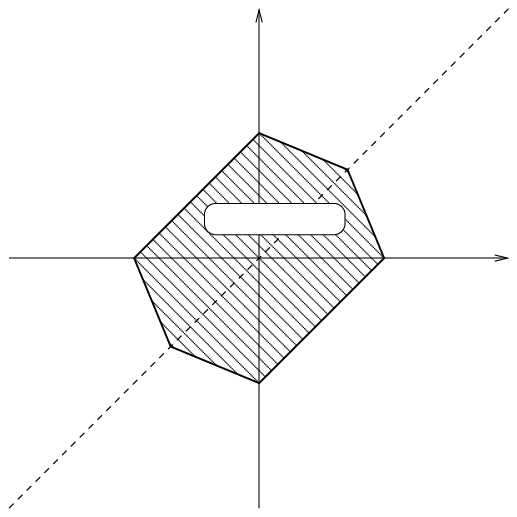}%
\end{picture}%
\setlength{\unitlength}{1973sp}%
\begingroup\makeatletter\ifx\SetFigFont\undefined%
\gdef\SetFigFont#1#2#3#4#5{%
  \reset@font\fontsize{#1}{#2pt}%
  \fontfamily{#3}\fontseries{#4}\fontshape{#5}%
  \selectfont}%
\fi\endgroup%
\begin{picture}(5100,4833)(-224,-3973)
\put(4876,-1636){\makebox(0,0)[lb]{\smash{{\SetFigFont{6}{7.2}{\familydefault}{\mddefault}{\updefault}{\color[rgb]{0,0,0}$\nu_1$}%
}}}}
\put(601,-361){\makebox(0,0)[lb]{\smash{{\SetFigFont{6}{7.2}{\familydefault}{\mddefault}{\updefault}{\color[rgb]{0,0,0}$|\nu_1|+|\nu_2|$}%
}}}}
\put(3826,-886){\makebox(0,0)[lb]{\smash{{\SetFigFont{6}{7.2}{\familydefault}{\mddefault}{\updefault}{\color[rgb]{0,0,0}$\sqrt{2}|\nu_2|+|\nu_2-\nu_1|$}%
}}}}
\put(3076,-3061){\makebox(0,0)[lb]{\smash{{\SetFigFont{6}{7.2}{\familydefault}{\mddefault}{\updefault}{\color[rgb]{0,0,0}$|\nu_1|+|\nu_2|$}%
}}}}
\put(2476,314){\makebox(0,0)[lb]{\smash{{\SetFigFont{6}{7.2}{\familydefault}{\mddefault}{\updefault}{\color[rgb]{0,0,0}$\sqrt{2}|\nu_1|+|\nu_2-\nu_1|$}%
}}}}
\put(526,-3661){\makebox(0,0)[lb]{\smash{{\SetFigFont{6}{7.2}{\familydefault}{\mddefault}{\updefault}{\color[rgb]{0,0,0}$\sqrt{2}|\nu_1|+|\nu_2-\nu_1|$}%
}}}}
\put(3751,-511){\rotatebox{45.0}{\makebox(0,0)[lb]{\smash{{\SetFigFont{6}{7.2}{\familydefault}{\mddefault}{\updefault}{\color[rgb]{0,0,0}$\nu_1=\nu_2$}%
}}}}}
\put(2416,704){\makebox(0,0)[lb]{\smash{{\SetFigFont{6}{7.2}{\familydefault}{\mddefault}{\updefault}{\color[rgb]{0,0,0}$\nu_2$}%
}}}}
\put(-224,-2041){\makebox(0,0)[lb]{\smash{{\SetFigFont{6}{7.2}{\familydefault}{\mddefault}{\updefault}{\color[rgb]{0,0,0}$\sqrt{2}|\nu_2|+|\nu_2-\nu_1|$}%
}}}}
\put(1891,-1261){\makebox(0,0)[lb]{\smash{{\SetFigFont{6}{7.2}{\familydefault}{\mddefault}{\updefault}{\color[rgb]{0,0,0}$F(p\cdot\Sigma)\le 1$}%
}}}}
\end{picture}%